\documentclass[11pt,twoside]{amsart}

\usepackage{amsfonts}

\title{On the use of K\"ulshammer type invariants in representation theory}
\address{Universit\'e de Picardie,\newline
D\'epartement de Math\'ematiques et UMR 6140 du CNRS, \newline 33 rue
St Leu,\newline F-80039 Amiens Cedex 1,\newline France}
\email{alexander.zimmermann@u-picardie.fr}
\author{Alexander Zimmermann}
\date{June 14, 2010}

\newtheorem{Theo1}{{Theorem}}
\newtheorem*{Theo2}{{Theorem}}
\newtheorem{Lemma1}{{Lemma}}[section]
\newtheorem{Def1}[Lemma1]{{Definition}}
\newtheorem{Prop1}[Lemma1]{{Proposition}}
\newtheorem{Claim1}[Lemma1]{{Claim}}
\newtheorem{Rem1}[Lemma1]{{Remark}}
\newtheorem{Cor1}[Lemma1]{{Corollary}}
\newtheorem{Ex1}[Lemma1]{{Example}}
\newtheorem*{Conj}{Conjecture}
\newtheorem*{Ack}{Acknowledgement}

\newenvironment{Lemma}{\begin{Lemma1}}{\end{Lemma1}}
\newenvironment{Def}{\begin{Def1}\em}{\end{Def1}}
\newenvironment{Prop}{\begin{Prop1}}{\end{Prop1}}
\newenvironment{Rem}{\begin{Rem1}\rm}{\end{Rem1}}
\newenvironment{Theorem}{\begin{Theo1}}{\end{Theo1}}

\newenvironment{Cor}{\begin{Cor1}}{\end{Cor1}}

\newenvironment{Example}{\begin{Ex1}\em}{\end{Ex1}}
\newenvironment{Acknowledgement}{\begin{Ack}\em}{\end{Ack}}

\newcommand{\lra}{\longrightarrow}

\newcommand{\ra}{\rightarrow}
\newcommand{\sdp}{\times\kern-.2em\vrule height1.1ex depth-.05ex}
\newcommand{\epi}{\lra \kern-.8em\ra}

\newcommand{\N}{{\mathbb N}}

\newcommand{\T}{{\mathbb T}}
\newcommand{\B}{{\mathbb B}}
\newcommand{\ul}{\underline}
\newcommand{\ol}{\overline}

\newcommand{\Z}{{\mathbb Z}}

\newcommand{\dickebox}{{\vrule height5pt width5pt depth0pt}}

 \normalbaselines
\begin{document}

\begin{abstract}
Since 2005 a new powerful invariant of an algebra emerged using earlier work of
Horv\'ath, H\'ethelyi, K\"ulshammer and Murray. The authors studied Morita
invariance of a sequence of ideals of the centre of a finite dimensional
algebra over a field of finite characteristic. It was shown that the sequence
of ideals is actually a derived invariant, and most recently a slightly modified
version of it an invariant under stable equivalences of Morita type. The invariant was
used in various contexts to distinguish derived and stable equivalence classes
of pairs of algebras in very subtle situations. Generalisations to non symmetric
algebras and to higher Hochschild (co-)homology was given. This article surveys
the results and gives some of the constructions in more detail.
\end{abstract}

\maketitle

\tableofcontents

\section*{Introduction}

Brauer studied representations of finite groups over fields of characteristic $p$
dividing the order of the group. In 1956 he showed \cite{Brauer}
amongst many other things that, if the field is algebraically closed,
the number of simple modules is equal to the number of conjugacy classes
of elements of $G$ of order prime to $p$. His method was rather
general already and
K\"ulshammer used these ideas 25 years later to define very
sophisticated invariants
for general symmetric algebras. More precisely K\"ulshammer
defined and studied in a series of four papers \cite{Ku1}
a sequence of ideals of the centre of a symmetric algebra.
In case of group rings over a finite group K\"ulshammer
studied lower defect
groups, proved Brauer's main theorems and many other results
known in modular
representation theory of finite groups
by these ideals and related invariants. K\"ulshammer's approach
was left untouched
until Murray \cite{Murray1,Murray2}
studied K\"ulshammer's approach in connection with the question of
the existence of real characters and of characters of defect $0$.
This was the starting point of the collaboration between Breuer,
H\'ethelyi, Horv\'ath,
K\"ulshammer and Murray \cite{HHKM} and \cite{BHHKM}
where the authors study the existence of odd diagonal entries in the
Cartan matrix
of a group algebra. Moreover, in \cite{HHKM} the Morita invariance of the
sequence of ideals is shown and the question of derived invariance is posed.

In \cite{Kuelsquest} completely different methods are used to show
that indeed at least for perfect base fields the sequence of ideals is
invariant under derived equivalence. Still the assumption that the algebras are
symmetric was needed.

After a lecture of the author in
September 2005 in Oberwolfach
Holm became interested in the sequence of ideals
and proposed several improvements and applications.
First, in joint work of the author with  Bessenrodt and Holm
\cite{nonsymmetric} using trivial extension algebras the
derived invariance, and also the very definition of the K\"ulshammer ideals
was extended to not necessarily symmetric algebras. Further, the K\"ulshammer
ideals were used in \cite{hz-tame}
to distinguish the derived equivalence class
of two algebras of dihedral type and of two pairs of algebras of semidihedral type
which were not seen to be not derived equivalent in the
classification of  Holm \cite{Holmhabil}.

During the September 2005 Oberwolfach lecture  Adem asked for a
generalisation of the K\"ulshammer ideal structure to higher Hochschild
(co-)homology. The question is far from trivial and was solved in \cite{Gerstenhaber}
where actually two approaches were taken, both of which are not exactly what
was asked for. The first approach uses Hochschild homology instead of cohomology
as in the K\"ulshammer ideal structure, the second approach uses the
Stasheff approach to the Gerstenhaber structure to get a non linear analogue.
The case of non symmetric algebras was answered in \cite{TAHochschild} again using
trivial extension algebras and the Hochschild homology approach.

The derived equivalence
classification of tame domestic weakly symmetric algebras
was given by Bocian, Holm and Skowro\'nski
for domestic algebras
\cite{BocianHolmSkowronski2004, BocianHolmSkowronski2007, weaklyprojdomestic} using
for the last remaining delicate questions K\"ulshammer ideals.
Similarly, also using K\"ulshammer ideals in parts,
Bia\l kowski, Holm and Skowro\'nski
\cite{BialkowskiHolmSkowronski2003a, BialkowskiHolmSkowronski2003b, weklyprojpolygrowth}
gave a derived equivalence classification
of tame algebras of polynomial growth up to some
difficult problems
concerning scalars in the relations of certain algebras, similar to the problem
solved in \cite{hz-tame}.

Then, most recently Bia\l kowski, Erdmann and Skowro\'nski classified selfinjective
algebras with the property that the third syzygy of every
simple module $S$ is again isomorphic to $S$. They obtain that these
algebras are all certain deformations of
preprojective algebras of a generalised Dynkin type. The deformations involve
parameters in the relations of the algebra. As was seen
\cite{weaklyprojdomestic,weklyprojpolygrowth,hz-tame}
K\"ulshammer ideals
are well suited for this kind of questions. Derived equivalence classes
of one family called of type $L$, defined in detail in
Example~\ref{ExampleFormIsNotSymmetric} and more generally in Section~\ref{BESdeformedpreprojective}
below, were largely
given in joint work with Holm \cite{preproj}.
Here we developed quite sophisticated methods to determine the
K\"ulshammer subspace structure at the beginning
which hold for a priori non symmetric algebras as well.
We display this method, even though strictly speaking
only small parts of it are really necessary in case of algebras of type $L$.
Nevertheless, the method works in general, is potentially very useful
and it seems reasonable to present it here.

During a lecture of the author in October 2007 at
Beijing Normal University the question of an
invariance under stable equivalences was posed.
In most recent results with Yuming Liu and Guodong Zhou a
generalisation of K\"ulshammer ideal theory was given for
stable categories and an
invariance was proved for stable equivalences of Morita type (\cite{LZZ}
and \cite{KLZ}). Most interestingly the result has strong links to the
Auslander-Reiten conjecture \cite[Conjecture 5, page 409]{ARS}, which
says that a stable equivalence should
preserve the number of simple non projective modules.
The result \cite{LZZ, KLZ} was used in joint work with Zhou \cite{Modulestructure}
to give a classification of algebras of
dihedral, of semidihedral and of quaternion type as defined by
Erdmann \cite{ErdmannLNM} up to stable equivalences of Morita type.
Moreover, in \cite{polygrowth}, in joint work with Zhou we prove
that the classification of weakly symmetric tame algebras of
polynomial growth
up to stable equivalence of Morita type coincides with the derived
equivalence classification of Bia\l kowski-Holm-Skowro\'nski.

In the present paper we survey these results and give at certain points quite
complete proofs for results
which seem to be very useful also in further contexts.
In Section~\ref{History} we trace some steps in the
origins, starting from Brauer and Reynolds.
Section~\ref{practicalquestions} reviews properties of selfinjective algebras
and develops tools to actually compute
K\"ulshammer ideals for quite complicated algebras. These tools were
developed during the past years for this purpose, but the origins are,
of course, very classical. Section~\ref{Kuelshammeridealtheory}
presents K\"ulshammer ideals as they were developed originally by K\"ulshammer and
as they were generalised later to non-symmetric algebras.
Section~\ref{Moritaderivedstablecategories} displays the
invariance of the various forms of K\"ulshammer ideals under Morita,
derived and stable equivalences.
Section~\ref{applisection} is devoted to various applications mentioned
above including a detailed outline of the
proof for the deformed preprojective
algebras of type $L$. Section~\ref{Hochschildsection} gives the above mentioned
two approaches
to the Hochschild (co-)homology generalisations of K\"ulshammer invariants.

\section{Historical facts and basic definitions: the origins by Brauer and Reynolds}

\label{History}\label{history}

Brauer developed in the 1950's the far reaching representation theory
of groups over fields of finite characteristic. In 1956 he showed in particular

\begin{Theorem}\label{Brauerstheorem} (Brauer \cite[Statement 3B]{Brauer})
If $K$ is an algebraically closed field of characteristic $p>0$ and
if $G$ is a finite group, then the number of simple $KG$-modules
up to isomorphism equals the number of conjugacy classes of $G$ of elements
of order prime to $p$.
\end{Theorem}

Of course, in the spirit of that time, Brauer did speak of irreducible characters
rather than of modules, but the result translates into modern terms as it is shown
above. The method of proof he used is somewhat indirect. He defines
for any $K$-algebra $A$ the space of commutators $[A,A]$, which is defined as the
$K$-vector space generated as vector space
by all possible expressions $ab-ba$ where $a,b\in A$.

Further he defines
$$TA:=\{a\in A\;|\;\exists n\in\N: a^{p^n}\in [A,A]\}$$
The first lemma is a little more general than stated in \cite{Brauer}, but
with identical proof. The generalisation comes from the fact that actually
one can consider more generally
$$T_nA:=\{a\in A\;|\; a^{p^n}\in [A,A]\}$$
for all $n\in\N$, a fact which is an observation due to K\"ulshammer.
We will need and study $T_nA$
later in more detail.

\begin{Lemma} (Brauer \cite[Statement 3A]{Brauer})\label{Brauerslemma}
Let $K$ be a field of characteristic $p>0$ and let $A$ be a $K$-algebra.
Then $T_nA$ is a $K$-subspace of $A$ satisfying
$TA=\bigcup_{n\in\N}T_nA$. If $A$ is finite dimensional, and $K$ is
a splitting field for $A$, then the number of
simple $A$-modules up to isomorphism equals the dimension of $A/TA$.
\end{Lemma}

The proof of this lemma is so simple that we may give it here in almost full detail.

\medskip

Proof (Brauer). Take $x,y\in T_nA$. Then develop $(x+y)^p$ and get a sum of
all possible words in $x$ and $y$ with $p$ factors, each occurring exactly once.
If $1<s<p$, then there are $n(s)$ of such words in which $x$ occurs $s$ times
and $y$ occurs $p-s$ times.
Take $N(s)$ the set of these words.
Then the cyclic group of order $p$ acts on this set by
cyclic permutation of the word:
$c\cdot (a_1a_2\dots a_{p-1}a_p):=(a_2a_3\dots a_pa_1)$
for a generator $c$ of the cyclic group and $a_i\in\{x,y\}$.
Hence $N(s)$ decomposes into orbits of length $p$, and the
difference of two elements in the same orbit is clearly
in $[A,A]$. Hence,
$$(x+y)^p-x^p-y^p\in[A,A].$$
Moreover
$$(xy-yx)^p+[A,A]=(xy)^p-(yx)^p+[A,A]=
x\left((yx)^{p-1}y\right)-\left((yx)^{p-1}y\right)x+[A,A]=[A,A]$$
and
$$(\lambda x)^p=\lambda^px^p$$
for all $x,y\in A$ and $\lambda\in K$ show that $T_nA$ is a $K$-subspace of $A$.

Traces of commutators of matrices are $0$.
Therefore, if $A=Mat_n(K)$, then
$[A,A]\subseteq \{M\in Mat_n(K)\;|\;trace(M)=0\}$.
On the other hand, using elementary matrices, one sees
that the inclusion actually is an equality. The space of
matrices with trace $0$ is of codimension $1$, and $Mat_n(K)$
has exactly one simple module up to isomorphism.
We denote by $rad(A)$ the Jacobson
radical of the algebra $A$. Hence, putting
$\ol A:=A/rad(A)$, one gets
$$\ol A/(T\ol A)=A/(TA)$$
as vector spaces since $rad(A)$ is nilpotent, and therefore
$rad(A)\subseteq T(A)$.
This shows the statement by Wedderburn's theorem. \dickebox

\medskip

Brauer's Theorem~\ref{Brauerstheorem} follows from the fact that for a
group ring $KG$
a basis is formed by the elements of $G$, and any element $g\in G$ admits a unique
(so-called the $p$-primary decomposition)
$g=g_p\cdot g_{p'}$ where $g_p$ is a $p$-element and $g_{p'}$ is
of order prime to $p$ commuting with $g_p$.
Hence
$$(g-g_{p'})^{p^n}+[KG,KG]=(g^{p^n}-g_{p'}^{p^n})+[KG,KG]=[KG,KG]$$
for a certain large $n$. Then,
$$hgh^{-1}-g=hgh^{-1}-gh^{-1}h=[h,gh^{-1}]\;\;\forall g,h\in G$$
and the rest is straightforward. \dickebox

\begin{Rem}\label{moreisshown}
Observe that actually much more is shown: $TA=rad(A)+[A,A]$.
\end{Rem}

Another concept due to Reynolds \cite{Reynolds} is closely linked.
Let $G$ be a finite group, let $K$ be an algebraically closed field of
characteristic $p>0$. For any $g\in G$ and $h\in G$ consider the $p$-primary decomposition
$g=g_p\cdot g_{p'}$ and $h=h_p\cdot h_{p'}$.
Then
$$S_h:=\{g\in G\;|\;\exists x\in G:x\cdot g_{p'}\cdot x^{-1}=h_{p'}\}$$
be the set of elements in $G$ whose $p'$-part is conjugate to the $p'$-part
of $h$.
Put $$C_h:=\sum_{g\in S_h}g$$
the sum of all these elements in $S_h$. Recall that the centre of $KG$
has a basis
consisting of all conjugacy class sums of elements of $G$.
Therefore $C_h\in Z(KG)$ for all $h\in G$.

\begin{Def} (Reynolds \cite[Theorem 1]{Reynolds})
The Reynolds ideal of $KG$ is the ideal of $Z(KG)$ generated as $K$-vector space by
the elements $C_h$, for $h\in G$.
\end{Def}

We get

\begin{Prop} (Reynolds \cite{Reynolds}, cf \cite[Theorem VI.4.6]{Feit})
The Reynolds ideal $R(KG)$ of $KG$ is the annihilator of $rad(KG)$ in $Z(KG)$.
\end{Prop}

\section{Selfinjective and symmetric algebras revisited}

\label{practicalquestions}

In order to be able to explain more deeply the relations between
Reynolds ideals, $T(KG)$ and related objects we need to
explain the structure of selfinjective and of symmetric algebras.
The theory is classical and originates in Nakayama's work
\cite{Nakayama} in the late 1930's.

Various approaches can be found in the literature, but as far as we know
the approach using Picard groups, which we will explain in
Section~\ref{Nakayamatwistedcentre} did not
appear elsewhere, though Yamagata~\cite{Yamagata} gives some related
thoughts.

Throughout this section we suppose for simplicity
that $K$ is a field. However, many
results stay true under weaker assumptions on $K$, sometimes
$K$ being a commutative ring would be sufficient.

\subsection{Basic definitions and properties}

Recall that for a $K$-algebra $A$ the space of linear forms $Hom_K(A,K)$
is an $A-A$-bimodule by
$$(afb)(x):=f(bxa)\;\;\forall a,b,x\in A\forall f\in Hom_K(A,K).$$
A group algebra is a symmetric algebra in the following sense.

\begin{Def}
Let $K$ be a field and let $A$ be a $K$-algebra. Then $A$ is
\begin{itemize}
\item symmetric if  $A\simeq Hom_K(A,K)$ as $A-A$-bimodules.
\item selfinjective if $A\simeq Hom_K(A,K)$ as $A$ left-modules.
\end{itemize}
\end{Def}

We shall derive some consequences.

Suppose $A$ is selfinjective and let $\varphi:A\lra Hom_K(A,K)$ be an isomorphism
of $A$ left-modules.
Then we may define a $K$-bilinear form
$$\langle \;,\;\rangle:A\times A\lra K$$
by
$$\langle a,b\rangle:=(\varphi(b))(a).$$
The fact that $\varphi$ is an isomorphism of vector spaces
is equivalent to the fact that $\langle \;,\;\rangle$ is non degenerate.

The fact that $\varphi$ is $A$-linear is equivalent to
$$\langle a,bc\rangle=(\varphi(bc))(a)=(b\varphi(c))(a)=\varphi(c)(ab)=\langle ab,c\rangle$$
for all $a,b,c$, where the linearity is used in the second equality. A bilinear form on an
algebra $A$ is called associative if $$\langle a,bc\rangle=\langle ab,c\rangle\mbox{
for all $a,b,c\in A$.}$$

Now $\varphi$ is an $A-A$-bimodule homomorphism
if and only if $\langle\;,\;\rangle$ is associative
(i.e. $\varphi$ is left $A$-linear) and moreover
$$\langle a,b\rangle=(\varphi(b))(a)=(\varphi(1)b)(a)=\varphi(1)(ba)=\langle ba,1\rangle=\langle b,a\rangle$$
and so $A$ is symmetric if and only if the associative
non degenerate form $\langle\;,\;\rangle$ may be chosen symmetric.

We summarise the statements in a (well known) proposition which gives an alternative
definition of selfinjective and symmetric algebras.

\begin{Prop}\label{formversusiso}\label{explicitselfinjective} Let $K$ be a field and
let $A$ be a finite dimensional $K$-algebra. Then we have the following statements.
\begin{itemize}
\item The algebra $A$ is selfinjective
if and only if there is a non degenerate associative bilinear form on $A$.
\item The algebra $A$ is symmetric if and only if there is a non degenerate
associative and symmetric bilinear form on $A$.
\end{itemize}
\end{Prop}

We should mention that the existence statement in Proposition~\ref{formversusiso}
is constructive: The bilinear form
is as explicit as is the isomorphism. That is, if one knows the explicit
isomorphism of $A$ to its dual by an explicit formula, then one knows the bilinear form
by an explicit formula, and vice versa.

The non degenerate associative symmetric bilinear form is called
symmetrising form for a symmetric algebra. In the remaining parts of Section~\ref{practicalquestions}
these ideas are developed further in particular with emphasis on
the question how to actually determine the bilinear form and associated questions.

For the moment we shall continue with Reynolds ideals and give the promised link.

In the following we frequently use for a symmetric algebra $A$ and subsets $S$ of $A$
the symbol $S^\perp$ to designate the
orthogonal space with respect to the symmetrising form of the algebra.

\begin{Prop} (K\"ulshammer~\cite{Ku2};
\cite[Part I Lemma A; Satz C; Satz D]{Ku1})\label{KAperp}
Let $K$ be a field and let $A$ be a finite dimensional symmetric $K$-algebra. Then $[A,A]^\perp=Z(A)$
and $soc(A)=rad(A)^\perp=Ann_A(rad(A))$. In particular $R(KG)=Z(KG)\cap soc(KG)$
for a finite group $G$.
\end{Prop}

Proof.
$$\langle ab-ba,c\rangle=\langle ab,c\rangle -\langle ba,c\rangle=
\langle a,bc\rangle -\langle c,ba\rangle=\langle bc,a\rangle
-\langle cb,a\rangle=\langle bc-cb,a\rangle$$
and hence $c\in [A,A]^\perp$ if and only if $\langle ab-ba,c\rangle=0$ for all $a,b$.
Therefore $c\in [A,A]^\perp$ if and only if $\langle bc-cb,a\rangle=0$ for all $a,b$.
In particular $c\in [A,A]^\perp$ if and only if $bc-cb\in A^\perp$. But $A^\perp=0$ since
the form is non degenerate. Hence $c\in [A,A]^\perp$ if and only if $bc=cb$
for all $b\in A$. This shows $[A,A]^\perp=Z(A)$.

$$\langle I,rad(A)\rangle=\langle 1,I\cdot rad(A)\rangle$$
and hence
$$I\subseteq rad(A)^\perp\Leftrightarrow I\cdot rad(A)=0\Leftrightarrow I\subseteq soc(A).$$
This finishes the proof.
\dickebox

\subsection{The Nakayama automorphism}\label{Nakayamatwistedcentre}

Selfinjective algebras come along with an automorphism, called
the Nakayama automorphism which will be explained now.

If $A$ is a selfinjective $K$-algebra, then
$A\simeq Hom_K(A,K)$ as an $A$ left-module. Hence, $Hom_K(A,K)$
is a free left $A$-module of rank $1$. Moreover,
$$End_A(\ _AHom_K(A,K))\simeq End_A(\ _AA)\simeq A$$
and so $Hom_K(A,K)$ is a progenerator over $A$ with endomorphism
ring isomorphic to $A$, whence inducing a Morita self-equivalence of $A$.
This implies that the isomorphism class of
$Hom_K(A,K)$ is in the Picard group $Pic_K(A)$ (cf e.g. \cite[Section 37]{MO}).
As is shown there there is a group homomorphism
$$\omega_0:Aut_K(A)\lra Pic_K(A)$$
given by $\omega_0(\alpha)=\ _1A_\alpha $. Here for any two
automorphisms $\alpha$ and $\beta$ of $A$ the $A-A$-bimodule
$_\alpha A_\beta$ denotes
$A$ as vector space, on which $a\in A$ acts by multiplication
by $\beta(a)$ on the right and by $\alpha(a)$ on the left.
To shorten the notation we abbreviate in this context the identity on $A$ by $1$.
One gets $ker(\omega_0)=Inn(A)$
the inner automorphisms of $A$ and hence
$$Out_K(A):=Aut_K(A)/Inn(A)$$
is a subgroup of $Pic_K(A)$. (Observe that the group of
inner automorphisms does not depend on $K$.)
The image of $\omega_0$ consists
of those isomorphism classes of invertible $A-A$-bimodules which are free
on the left (cf. \cite[(37.16) Theorem]{MO}). Observe that
\begin{eqnarray*}
\ _{\alpha^{-1}} A_1&\stackrel{}{\lra}&\ _1A_\alpha\\
a&\mapsto&\alpha(a)
\end{eqnarray*}
is an $A-A$-bimodule homomorphism. Hence $\omega_0$ may also be
defined by twisting the action on the left.

Now, as $Hom_K(A,K)$ is free of rank $1$
as left-module one gets that $Hom_K(A,K)$ is in the image of
$Out_K(A)$ in $Pic_K(A)$ and
therefore there is an automorphism $\nu\in Aut_K(A)$ so that
$$Hom_K(A,K)\simeq \ _1A_\nu $$
as $A$-$A$-bimodules. The automorphism
$\nu$ is unique up to an inner automorphism.

\begin{Def}
Let $A$ be a selfinjective $K$-algebra. Then there is an automorphism
$\nu$ of $A$ so that $Hom_K(A,K)\simeq \ _1A_\nu$ as $A$-$A$-bimodules.
This automorphism is unique up to inner automorphisms and is called the
{\em Nakayama automorphism}.
\end{Def}

\begin{Rem}
In Nakayama's original approach $_\nu A_1$ is used instead of $_1A_\nu$ and so
the Nakayama automorphism in Nakayama's work corresponds to the inverse of what
we define here.
\end{Rem}

In principle, the definition of selfinjectiveness uses an isomorphism
of the $A$-left-module of linear forms on $A$ with the regular $A$-module.
One could use the right module structure as well. We get the well-known result
that left selfinjective is equivalent to right selfinjective.

\begin{Cor}
$_AA\simeq\ _AHom_K(A,K)\Leftrightarrow A_A\simeq Hom_K(A,K)_A$.
\end{Cor}

To prove the Corollary one just needs to see that the
isomorphism as left-modules implies the following bimodule isomorphisms.
$$Hom_K(A,K)\simeq \ _1A_{\nu}\simeq \ _{\nu^{-1}}A_1$$
and so $Hom_K(A,K)\simeq A$ as $A$ right-modules.

\medskip

Now, for a selfinjective $K$-algebra $A$, given a simple $A$-module $S$, then
$_1A_\nu \otimes_AS\simeq Hom_K(A,K)\otimes_AS$ is again a simple $A$-module.

\begin{Def}
Let $K$ be a field and let $A$ be a finite dimensional $K$-algebra. Then
$A$ is {\em weakly symmetric} if $A$ is selfinjective and
$Hom_K(A,K)\otimes_AS\simeq S$ for all simple $A$-modules $S$.
\end{Def}

This definition will be important in Section~\ref{polygrowthsection}.

\subsection{The Nakayama twisted centre}
\label{Nakayamatwistedcentresection}

Let $K$ be a field and let $A$ be a finite dimensional $K$-algebra. For an explicitly given
algebra, say as quiver with relations,
it is not very hard to write down many commutators. This gives an upper bound
for the dimension of $A/[A,A]$. However, to prove that
the commutators found really generate $[A,A]$ is quite difficult
in general. The method is to interpret $A/[A,A]$ as
different space, in which it is easier to find many linearly independent
elements. This then gives a lower bound for the dimension of $A/[A,A]$.
If the lower and the upper bound coincide, then one has proved that the
commutators found actually generate the whole space $[A,A]$.

Let $A$ be a selfinjective $K$-algebra. Given a ring $R$ and a right
$R$-module $M$ and a left $R$-module $N$, the very definition of
the tensor product $M\otimes_RN$ as free abelian on symbols $m\otimes n$
with relations $m\otimes rn=mr\otimes n$ and additivity in each variable
gives an isomorphism $A/[A,A]\simeq A\otimes_{A\otimes_KA^{op}}A$.
An alternative way to see this is by the bar resolution of Hochschild homology.

Now
\begin{eqnarray*}
Hom_K(A/[A,A],K)
&\simeq & Hom_K(A\otimes_{A\otimes_KA^{op}}A,K)\\
&\simeq& Hom_{A\otimes_KA^{op}}(A, Hom_K(A,K))\\
&\simeq &Hom_{A\otimes_KA^{op}}(A,\ _1A_\nu )
\end{eqnarray*}
which gives
\begin{eqnarray*}
Hom_K(A/[A,A],K)
&\simeq &\{a\in A\;|\;b\cdot a=a\cdot \nu(b)\;\forall b\in A\}
\end{eqnarray*}
where the isomorphism is given by sending a homomorphism to
the image of $1\in A$ which will satisfy the equation by
the property of the homomorphism being $A\otimes_KA^{op}$-linear.

\begin{Def} (Holm and Zimmermann \cite[Definition 2.2]{preproj})
Let $K$ be a field and let $A$ be a selfinjective $K$-algebra with
Nakayama automorphism $\nu$. Then the {\em Nakayama twisted centre} is
defined to be
$$Z_\nu(A):=\{a\in A\;|\;b\cdot a=a\cdot\nu(b)\;\forall b\in A\}.$$
\end{Def}

\begin{Rem}
The definition works for $K$ a commutative ring as well.

The automorphism $\nu$ is unique only up to an inner automorphism.
If $\nu$ differs from $\nu'$ by an inner automorphism,
let $\nu(b)=u\cdot \nu'(b)\cdot u^{-1}$ for all $b\in A$ and
some unit $u$ of $A$. Then
\begin{eqnarray*}
Z_\nu(A)&=&\{a\in A\;|\;b\cdot a=a\cdot\nu (b)\;\forall b\in A\}\\
&=&\{a\in A\;|\;b\cdot a=a\cdot u\cdot \nu'(b)\cdot u^{-1}\;\forall b\in A\}\\
&=&\{a\in A\;|\;b\cdot (a\cdot u)=(a\cdot u)\cdot \nu'(b)\;\forall b\in A\}\\
&=&\{a\in A\;|a\cdot u\in Z_{\nu'}(A)\}\\
&=&Z_{\nu'}(A)\cdot u^{-1}
\end{eqnarray*}
\end{Rem}

\begin{Rem}
In general the Nakayama twisted centre will not be a ring: if $a,b\in Z_\nu(A)$,
then $$b(a_1a_2)=(ba_1)a_2=(a_1\nu(b))a_2=a_1(\nu(b)a_2)=a_1a_2\nu^2(b)$$
and $\nu^2=\nu$ is equivalent to $\nu=id$.
Nevertheless, if $z\in Z(A)$ and $a\in Z_\nu(A)$ then
$$b\cdot za=zba=za\cdot \nu(b)$$
and $za\in Z_\nu(A)$. Hence $Z_\nu(A)$
is a $Z(A)$-submodule of $A$. The module structure does not depend on the
chosen Nakayama automorphism, up to isomorphism of $Z(A)$-modules.
\end{Rem}

\begin{Rem}
In case one follows Nakayama's original definition of a Nakayama
automorphism we need to replace $\nu$ by $\nu^{-1}$ and
hence there the Nakayama twisted centre would consist of elements $a$
satisfying $\nu(b)\cdot a=a\cdot b$.
\end{Rem}

We summarise our results in the following

\begin{Prop}\label{propertiesofNakayamatwistedcentre}
(Holm and Zimmermann \cite[Lemma 2.4]{preproj})
If $A$ is a selfinjective $K$-algebra over a field $K$ with Nakayama automorphism $\nu$.
Then $Hom_K(A/[A,A],K)\simeq Z_\nu(A)$ as $Z(A)$-modules. \dickebox
\end{Prop}

Again, the proposition holds as well for $K$ being a commutative ring.

\subsection{How to get the Nakayama automorphism explicitly}
\label{Nakayamainpractice}

Let $K$ be a field and let $A$ be a finite dimensional selfinjective $K$-algebra.

In order to compute the Nakayama automorphism $\nu$ we need to find an explicit
isomorphism $A\lra Hom_K(A,K)$ as $A$-modules.

\begin{Prop}\label{thenuformula} (Holm and Zimmermann \cite[Lemma 2.7]{preproj})
Let $K$ be a field and let $A$ be a finite dimensional selfinjective $K$-algebra
with associated bilinear form $\langle\;,\;\rangle$.
Then the Nakayama automorphism $\nu$ of $A$,
satisfies $\langle a,b\rangle=\langle b,\nu(a)\rangle$ for all $a,b\in A$,
and any automorphism satisfying this formula is a Nakayama automorphism.
\end{Prop}

\begin{Rem}
If one would use Nakayama's original definition the form would
satisfy $\langle \nu(a),b\rangle=\langle b,a\rangle$ for all $a,b\in A$.
\end{Rem}

Proof of Proposition~\ref{thenuformula}.
There is a non-degenerate associative bilinear
form on $A$, which induces an isomorphism between $A$ and the linear forms
on $A$ as $A$-modules by Proposition~\ref{formversusiso}.
The isomorphism gives an isomorphism of $A$-$A$-bimodules of
$ _1A_\nu$ and $Hom_K(A,K)$
by
\begin{eqnarray*}
_1A_\nu&\stackrel{\varphi}{\lra}& Hom_K(A,K)\\
a&\mapsto&\langle -,a\rangle=\varphi(a).
\end{eqnarray*}
Therefore, $\varphi(a)=\varphi(1\cdot \nu^{-1}(a))=\varphi(1)\cdot \nu^{-1}(a)$
and $\varphi(a)=\varphi(a\cdot 1)=a\cdot\varphi(1)$.
Since for $f\in Hom_K(A,K)$ one has $(fa)(b)=f(ab)$ and $(af)(b)=f(ba)$
for all $a,b\in A$, one gets
\begin{eqnarray*}
\langle b,a\rangle&=&(\varphi(a))(b)
=(\varphi(1)\cdot\nu^{-1}(a))(b) = \varphi(1)(\nu^{-1}(a)b)\\
&=&(b\cdot\varphi(1))(\nu^{-1}(a)) = \varphi(b)(\nu^{-1}(a))
 = \langle \nu^{-1}(a),b\rangle.
\end{eqnarray*}
Now, putting $a:=\nu(a')$ one gets
$$\langle b,\nu(a')\rangle=\langle a',b\rangle$$
for all $a',b\in A$.
Hence, the Nakayama automorphism has the above property.
Conversely, if an automorphism $\nu$ satisfies
$\langle a,b\rangle=\langle b,\nu(a)\rangle$ for all $a,b\in A$,
then the mapping $A\lra Hom_K(A,K)$ given by $a\mapsto\langle -,a\rangle$ gives an isomorphism
of $A$ and $Hom_K(A,K)$ as $A$-modules, inducing the element $_1A_\nu$ in the
Picard group of $A$.
\dickebox

\subsection{Practical questions for algebras given by quivers and relations}

In \cite{hz-tame} the following very useful result was proved for weakly
symmetric algebras. However the statement holds in a more general form.

\begin{Prop}
\label{prop:form}
Let $K$ be a field and let $A=KQ/I$ be a selfinjective algebra given by the quiver $Q$ and
ideal of relations $I$, and fix
a $k$-basis ${\mathcal B}$ of $A$ consisting of pairwise
distinct non-zero paths of the quiver $Q$. Assume that
${\mathcal B}$ contains a basis of the socle $soc(A)$ of $A$.
Define a $K$-linear mapping $\psi$ on the basis elements by
$$
\psi(b)=\left\{
\begin{array}{ll} 1 & \mbox{if $b\in soc(A)$} \\
                  0 & \mbox{otherwise}
\end{array} \right.
$$
for $b\in {\mathcal B}$.
Then an associative non-degenerate $K$-bilinear
form $\langle-,-\rangle$ for $A$ is given by
$\langle x,y\rangle := \psi(xy).$
\end{Prop}

\begin{Rem}
In case $A$ is weakly symmetric Proposition~\ref{prop:form}
was proved in \cite{hz-tame}. The assumption that $A$ is weakly symmetric was
used in \cite{hz-tame} only to prove the non degeneracy
of the form. For the reader's convenience we include a complete
proof.
\end{Rem}

Proof.
By definition, since $A$ is an associative algebra,
$\psi$ is associative on basis elements,
hence is associative on all of $A$.

Let $\nu$ be a Nakayama automorphism of $A$.
We observe now that $\psi(x\cdot \nu(e))=\psi(e\cdot x)$ for all $x\in A$ and all
primitive idempotents $e\in A$. Indeed, since $\psi$
is linear, we need to show this only on the elements
in $\mathcal B$. Let $b\in\mathcal B$. If $b$ is a path not in the
socle of $A$, then $b\nu(e)$ and $eb$ are either zero or not contained
in the socle either,
and hence $0=\psi(b)=\psi(b\nu(e))=\psi(eb)$.
If $b\in{\mathcal B}$ is in the socle of $A$, then $b=e_bb=b\nu(e_b)$ for exactly
one primitive idempotent $e_b$ and $e'b=b\nu(e')=0$ for each primitive idempotent
$e'\neq e_b$. Therefore, $\psi(e'b)=\psi(b\nu(e'))=0$ and
$\psi(e_bb)=\psi(b)=\psi(b\nu(e_b))$.

It remains to show that the map $(x,y)\mapsto \psi(xy)$ is non-degenerate.
Suppose we had $x\in A\setminus\{0\}$ so that $\psi(xy)=0$ for all $y\in A$.
In particular for each primitive idempotent $e_i$ of $A$ we get
$\psi(e_ixy)=\psi(xy\nu(e_i))=0$
for all $y\in A$.
Hence we may suppose that
$x\in e_iA$
for some primitive idempotent $e_i\in A$.

Now, $xA$ is a right $A$-module. Choose a simple submodule $S$ of $xA$ and
$s\in S\setminus\{0\}$. Then, since $s\in S\leq xA$ there is a $y\in A$ so that $s=xy$.
Since $S\leq xA\leq A$, and since $S$ is simple, $s\in soc(A)\setminus\{0\}$.
Moreover, since $x\in e_iA$, also $s=e_is$, i.e. $s$ is in the (1-dimensional)
socle of the projective indecomposable module $e_iA$. So, up to a non zero scalar factor,
$s$ is a path contained in the basis $\mathcal{B}$ (recall that
by assumption $\mathcal{B}$ contains a basis of the socle).
This implies that
$$\psi(xy)=\psi(s)=\psi(e_is)\neq 0,$$
contradicting the choice of $x$, and hence proving non-degeneracy. \dickebox

\begin{Rem}
It should be noted that the form depends on the chosen basis of the algebra.
Indeed, take $A=K[X]/X^2$. The socle is one-dimensional, and take a basis $\{X\}$. Then
one may complete with an element $1+\mu X$ for $\mu\in K$ to a basis of $A$. Hence
$1+X=(1+\mu X)+(1-\mu) X$ and
we get $\langle 1+X,1\rangle=1-\mu$ which depends heavily on $\mu$.
\end{Rem}

\begin{Example}\label{twistedcentreexample}
This example was communicated to me by Guodong Zhou
in January 2010 during a visit in Paderborn.
Let $K$ be a field and let $A_q=K\langle X,Y\rangle/(X^2,Y^2,XY-qYX)$ for $q\neq 0$.
Then $A_q$ is always selfinjective and $A_q$ is
symmetric if and only if $q=1$. Now, $\nu(X)=qX$ and $\nu(Y)=q^{-1}Y$ defines a Nakayama
automorphism. Indeed if
we use the $K$-basis $\{1,X,Y,XY\}$ for $A_q$ we get
$$1=\langle X,Y\rangle=\langle Y,\nu(X)\rangle=\langle Y,qX\rangle=1$$
$$1=\langle Y,X\rangle=\langle X,\nu(Y)\rangle=\langle X,q^{-1}Y\rangle=1$$
and likewise for the other basis elements.

Now, if $q\neq 1$, we get $$Z(A_q)=K\cdot 1+K\cdot XY$$ whereas
$$Z_\nu(A_q)=K\cdot X+K\cdot Y+K\cdot XY.$$
Hence the $Z(A_q)$-module $Z_\nu(A_q)$ is isomorphic to $\left(Z(A_q)/rad(Z(A_q))\right)^3$.

If $A$ is a symmetric algebra, and $\nu$ is chosen to be inner,
then $Z_\nu(A)$ is a rank $1$ free $Z(A)$-module.

Observe that the order of the above automorphism is the multiplicative order of
$q$ in $K$. Hence for big fields $K$ it is possible to create algebras with Nakayama
automorphisms of any given order, even infinite order.
\end{Example}

\begin{Example}\label{ExampleFormIsNotSymmetric}
Proposition~\ref{prop:form} is stated in \cite{weaklyprojdomestic} and in \cite{weklyprojpolygrowth} with an additional conclusion. Namely it is stated
there that the form $\langle\;,\;\rangle$ is symmetric in case the algebra is
symmetric. This is not true in general. A counterexample was given during the collaboration on \cite{preproj}.
Namely, the deformed preprojective algebra in the sense of Bia\l kowski,
Erdmann and Skowro\'nski \cite{BES}
of type $L_n$ for $n\geq 3$ gives an example (see Remark~\ref{BESmistakes} below).
This algebra is defined by the quiver
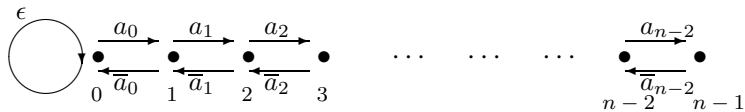
\unitlength1cm
\begin{center}
\begin{picture}(10,2)
\put(1,1){$\bullet$}
\put(1,0.5){\scriptsize $0$}
\put(2,1){$\bullet$}
\put(2,0.5){\scriptsize $1$}
\put(3,1){$\bullet$}
\put(3,0.5){\scriptsize $2$}
\put(4,1){$\bullet$}
\put(4,0.5){\scriptsize $3$}
\put(8,1){$\bullet$}
\put(7.8,0.4){\scriptsize $n-2$}
\put(9,1){$\bullet$}
\put(9,0.4){\scriptsize $n-1$}
\put(5,1){$\cdots$}
\put(6,1){$\cdots$}
\put(7,1){$\cdots$}
\put(1.1,1.3){\vector(1,0){.8}}
\put(2.1,1.3){\vector(1,0){.8}}
\put(3.1,1.3){\vector(1,0){.8}}
\put(8.1,1.3){\vector(1,0){.8}}
\put(8.9,.9){\vector(-1,0){.8}}
\put(3.9,.9){\vector(-1,0){.8}}
\put(2.9,.9){\vector(-1,0){.8}}
\put(1.9,.9){\vector(-1,0){.8}}
\put(0.4,1.1){\circle{1}}
\put(.88,1.1){\vector(0,-1){.1}}
\put(0,1.6){\small $\epsilon$}
\put(1.3,1.4){\small $a_0$}
\put(2.3,1.4){\small $a_1$}
\put(3.3,1.4){\small $a_2$}
\put(8.3,1.4){\small $a_{n-2}$}
\put(1.3,.65){\small $\ol a_0$}
\put(2.3,.65){\small $\ol a_1$}
\put(3.3,.65){\small $\ol a_2$}
\put(8.3,.65){\small $\ol a_{n-2}$}
\end{picture}
\end{center}
subject to the following relations
$$a_i\overline a_i+\overline a_{i-1}a_{i-1}=0\mbox{ for all }
i\in\{1,\dots ,n-2\}~,
$$
$$\overline a_{n-2}a_{n-2} =0~,~~~~~
\epsilon^{2n}=0~,~~~~~
\epsilon^2+a_0\overline a_0+\epsilon^{3}p(\epsilon)=0
$$
for a polynomial $p(X)\in K[X]$.
These algebras are the {\em deformed preprojective algebras of
type $L_n$}, in the sense of Bia\l kowski, Erdmann and Skowro\'{n}ski
\cite{BES}.

For the special case $p(X)=X^{2j}$ for $j\in\N$ and
for abbreviation we call this algebra by $L_n^j$ and assume that $K$ is of
characteristic $2$. Here we just give an example where the bilinear form
Proposition~\ref{prop:form} does not yield a symmetric bilinear form. We will deal with the
general case in a later Section~\ref{BESdeformedpreprojective}.

In order to be able to apply Proposition~\ref{prop:form} we need to fix a basis of
the socle of $L_n^j$. The fact that the elements below
is indeed a basis of the algebra is shown in \cite{preproj}
and the basis displayed in Proposition~\ref{abasis} can easily be transformed into the basis
below. Most recently in a completely independent approach
Andreu \cite{AndreuEstefania} shows that the basis below is indeed a basis.

For our purpose it seems to be most natural to take as $K$-basis of the socle the set
$$\{\epsilon^{2n-1},
\ol a_{i-1}\ol a_{i-2}\dots\ol a_0\epsilon^{2n-3-2i}
a_0a_1\dots a_{i-1}\;|\;i\in\{1,2\dots,n-2\}\}.$$
Complete the elements
$$\ol a_{i-1}\ol a_{i-2}\dots\ol a_0\epsilon^{2n-3-2i}
a_0a_1\dots a_{i-1}$$
of $e_iL_n^je_i$, for $i\geq 1$, to a basis of $e_iL_n^je_i$ by the elements
$$\ol a_{i-1}\ol a_{i-2}\dots\ol a_0\epsilon^{\ell}
a_0a_1\dots a_{i-1}$$
for $\ell\leq 2n-4-3i$ and
$$a_ia_{i+1}a_{i+2}\dots a_{j}\ol a_{j}\ol a_{j-1}\dots\ol a_{i+2}\ol a_{i+1}\ol a_i$$
for $i+1\leq j\leq n-2$.
A basis of $e_0L_n^je_0$ is given by $\epsilon^\ell$, for $0\leq \ell\leq 2n-1$.

Now we verify
\begin{eqnarray*}
\langle \ol a_0\epsilon^m,a_0\rangle&=&
\left\{\begin{array}{ll}1&\mbox{ if }m=2n-3\\
0&\mbox{ if }m\neq 2n-3\end{array}\right.\\
\langle a_0,\ol a_0\epsilon^m\rangle&=&
\begin{cases}1&\mbox{ if }m=2n-3\mbox{ or }m=2n-4-2j\\
 0&\mbox{ else}
\end{cases}
\end{eqnarray*}
Hence the bilinear form from Proposition~\ref{prop:form} is not symmetric.
However, the algebras $L_n^j$ are symmetric (cf Proposition~\ref{Lnissymmetric} below).
\end{Example}

\section{K\"ulshammer, the new idea in the 1980s}
\label{Kuelshammeridealtheory}

We come back to Brauer's proof displayed in Section~\ref{history}, Reynolds discoveries and present the
original approach which was introduced by K\"ulshammer to improve and unify these
earlier approaches. Moreover we explain the generalisation to non symmetric
algebras.

\subsection{K\"ulshammer's original construction for symmetric algebras}

Recall that for a symmetric $K$-algebra $A$ one defined
$T_n(A):=\{a\in A|\;a^{p^n}\in [A,A]\}$ and $T_n^\perp(A)$ is the orthogonal space
with respect to the symmetrising form.

\begin{Def} (K\"ulshammer~\cite[Part IV]{Ku1})
Let $A$ be a symmetric $K$-algebra. Then the ideal $T_n(A)^\perp$ of $Z(A)$ is
the $n$-th K\"ulshammer ideal.
\end{Def}

\begin{Rem}
\begin{itemize}
\item Remark~\ref{moreisshown} gives
that $T(A)=[A,A]+rad(A)$.  Hence
$$T(A)^\perp=rad(A)\cap soc(A)=:R(A).$$
\item
K\"ulshammer calls the ideals $T_n(A)^\perp$ the generalised Reynolds ideals.
\item
Since $T_n(A)\subseteq T_{n+1}(A)$
we get $T_n(A)^\perp\supseteq T_{n+1}(A)^\perp$
and the set of K\"ulshammer ideals is a decreasing
sequence of ideals of the centre
with first term being the centre and last term being the Reynolds ideal $R(A)$.
\end{itemize}
\end{Rem}

K\"ulshammer obtained in \cite{Ku1} many properties of a group in terms of properties of
the sequence of ideals $T_n(KG)^\perp$. Already in the first discussions in \cite{Ku1}
the setup was completely general and the definitions were given for a symmetric algebra
over a field of finite characteristic in general. Nevertheless, the applications in
focus in K\"ulshammer's discussions mainly have been group algebras and representations
of finite groups.

The main tool for technical proofs is still the symmetrising form.
Since we have by Proposition~\ref{KAperp} that $[A,A]^\perp=Z(A)$,
the restriction of the symmetrising form $\langle\;,\;\rangle$ to $Z(A)$
on the left argument
induces a non degenerate form, also denoted by $\langle\;,\;\rangle$,
$$\langle\;,\;\rangle:Z(A)\times A/[A,A]\lra K.$$
As we have seen in Lemma~\ref{Brauerslemma} the mapping
\begin{eqnarray*}
A/[A,A]&\stackrel{\mu}{\lra}&A/[A,A]\\
a+[A,A]&\mapsto&a^p+[A,A]
\end{eqnarray*}
is additive and semilinear, i.e. linear if one applies in addition a twist with the
Frobenius automorphism of the field.
If $V$ and $W$ are finite dimensional vector spaces over a field $K$
and if $\langle.,.\rangle$ is a non degenerate bilinear pairing
$V\times W\lra K$. Then any endomorphism $\varphi$ of $W$ has a unique left adjoint $\varphi^*\in End_K(V)$
satisfying $\langle v,\varphi(w)\rangle=\langle \varphi^*(v),w\rangle$ for all $v\in V$ and $w\in W$.

Now, this fact holds for semi-linear maps as
well as for linear maps, and the map $\mu$ has a left adjoint
$$\zeta:Z(A)\lra Z(A).$$
Now, $a\in T_n(A)\Leftrightarrow a\in ker(\mu^n)$
gives

\begin{Lemma}\label{zeta}
$z\in T_n(A)^\perp\Leftrightarrow z\in im(\zeta^n).$
\end{Lemma}
This characterisation will be the main tool for most of the
abstract statements later.

\subsection{Extending to general algebras: trivial extension algebras and K\"ulshammer theory}

Up to now, in order to establish a K\"ulshammer ideal theory it was necessary
already for the very definition to work over symmetric algebras. There is a
method to circumvent this difficulty.

Let $A$ be any finite dimensional $K$-algebra.
Then, as already mentioned, $Hom_K(A,K)$ is an $A-A$-bimodule by
$$(afb)(c):=f(bca)\;\;\forall a,b,c\in A\mbox{ and }f\in Hom_K(A,K).$$

Recall the construction of the trivial extension algebra, which is well-known
and very useful in the representation theory of associative algebras.
We may form $\T A:=Hom_K(A,K)\times A$, which is naturally a $K$-vector space.
We may define an algebra structure on this space.
$$(f,a)\cdot (g,b):=(ag+fb,ab)\;\;\forall\;a,b\in A,f,g\in Hom_K(A,K)\;.$$
It is a tedious but straightforward computation to verify that
$\T A$ is a $K$-algebra by this multiplication. Moreover
the projection to the second component is an algebra homomorphism $\T A\lra A$
with kernel $Hom_K(A,K)$ being an ideal with square $0$.

\begin{Def}
For any finite dimensional $K$-algebra $A$ the algebra $\T A$ is called the
trivial extension algebra.
\end{Def}

The property that is most interesting for our purposes is that $\T A$ is
a symmetric algebra, whatever may be the structure of the algebra $A$.
Indeed,
$$\langle (f,a),(g,b)\rangle:= g(a)+f(b)\;\;\forall (f,a),(g,b)\in\T A$$
is a symmetric associative non degenerate bilinear form on $\T A$.
This fact can be found in e.g. \cite[Section 3]{nonsymmetric}.
Proposition~\ref{formversusiso} shows then the statement.

In \cite{nonsymmetric} Bessenrodt, Holm and the author compute the
K\"ulshammer ideals of $\T A$.
We denote
$$Ann_{Hom_K(A,K)}(I):=\{f\in Hom_K(A,K)\;|\;f(I)=\{0\}\}$$
for any subset $I\subseteq A$. With this notation we showed

\begin{Prop} (Bessenrodt, Holm, Zimmermann \cite[Theorem 4.1]{nonsymmetric})
\label{prop-center}
Let $A$ be a finite-dimensional algebra
over a field of characteristic $p>0$, and let
${\T}A$ be its trivial extension.
\begin{enumerate}
\item[{(1)}] We have $T_0({\T}A)^\perp=Z({\T}A)=
Ann_{Hom_K(A,K)}([A,A])\times Z(A) $.
\item[{(2)}] For all $n\geq 1$ one has
$T_n({\T}A)^{\perp} = Ann_{Hom_K(A,K)}(T_nA) \times 0.$
\end{enumerate}
\end{Prop}

This result, though not difficult to prove,
is most remarkable, since for symmetric algebras we may
use the symmetrising form to transport $T_n(A)$ via orthogonality to
an ideal $T_n(A)^\perp$ of the centre of $A$. The ideal structure allows to
consider many invariants from commutative algebra attached to this
ideal $T_n(A)^\perp$. If $A$ is not symmetric, this is not easily possible. Hence
it is surprising that enlarging $A$ to $\T A$, the space $T_n(\T A)^\perp$
keeps the trace of $T_n(A)^\perp$ faithfully in the sense that $T_n(A)^\perp$ can be fully recovered
by $T_n(\T A)^\perp$.

If $A$ is already symmetric,
the isomorphism $\T A\lra Hom_K(\T A,K)$
takes $T_n(A)^\perp$ to $Ann_{Hom_K(A,K)}(T_n(A))$.
This fact is not hard to see. Indeed, if $A$ is symmetric,
$Hom_K(A,K)$ is the space of all $\langle a,-\rangle$ for
$a\in A$. Now, a linear form being in $Ann_{Hom_K(A,K)}(T_n(A))$
is equivalent to $\langle a,T_n(A)\rangle=0$ which in turn is equivalent
to $a\in T_n(A)^\perp$.

\section{Morita, derived and stable invariance}

\label{Moritaderivedstablecategories}

\subsection{Morita invariance}

In \cite{HHKM} H\'eth\'elyi, Horv\'ath, K\"ulshammer and Murray
studied amongst other questions the invariance of K\"ulshammer ideals
$T_n(A)^\perp$ under Morita equivalence. Recall that if
$$M\otimes_A-:A-mod\lra B-mod$$ is an equivalence, then
for any $z\in Z(A)$ there is a unique $\varphi_M(z)\in Z(B)$ so that
$m\cdot z=\varphi_M(z)\cdot m$ for all $m\in M$. Then,
$\varphi_M:Z(A)\lra Z(B)$ is an isomorphism of algebras.

\begin{Prop}(H\'ethelyi, Horv\'ath, K\"ulshammer and Murray \cite[Corollary 5.3]{HHKM})
Let $K$ be a perfect field of characteristic $p>0$ and let $A$ and $B$ be
finite dimensional $K$-algebras. If $$M\otimes_A-:A-mod\lra B-mod$$
is a Morita equivalence, then
$\varphi_M(T_n(A)^\perp)=T_n(B)^\perp$
for all $n\in\N$.
\end{Prop}

The authors show that
the mapping $\zeta$ of Lemma~\ref{zeta} behaves well with respect to multiplication
by idempotents. Using this statement it is possible to reduce to basic algebras, and to
use then that two
Morita equivalent basic algebras are isomorphic.

The existence of basic algebras, i.e. an up to isomorphism unique minimal
algebra
which is Morita equivalent to the given algebra, is very specific
for Morita equivalences.
Such a concept does not exist for weaker equivalences such as derived
equivalences
or stable equivalences of Morita type.

\subsection{Derived invariance}

\cite[Question 5.4]{HHKM} asked if K\"ulshammer ideals are also invariant under
derived equivalences. The method used for the Morita invariance does not apply since
as said before a concept of ''derived basic'' algebras do not exist.
Nevertheless, an equivalence between derived categories of finite
dimensional algebras imply the existence of an isomorphism of the centres.

More precisely, denote by $D^b(A)$ the derived category of bounded
complexes of finitely
generated $A$-modules with non zero homology in only finitely many degrees.

The main tool is the following result.

\begin{Theorem} (Rickard \cite[Theorem 3.3]{Ri1}, Keller \cite[Section 8]{Keller};
cf also e.g. \cite[Theorem 6.2.8]{derbuch})
Let $K$ be a field, let $A$ and $B$ be finite dimensional
$K$-algebras and suppose that
$D^b(A)\simeq D^b(B)$ is an equivalence of triangulated categories.
Then there is a complex $X\in D^b(B\otimes_KA^{op})$ which is formed
by modules which
are projective as $B$-modules and projective as $A$-modules, so that
$$X\otimes_A^{\mathbb L}-:D^b(A)\lra D^b(B)$$
is an equivalence of triangulated categories.
\end{Theorem}

A complex $X$ as in the theorem is called two-sided tilting complex. It is
unknown if every equivalence between derived categories is of the form
$$X\otimes_A^{\mathbb L}-:D^b(A)\lra D^b(B)$$
This is shown to hold on the level of objects, but it is not known if there may be
an exotic equivalence behaving differently on morphisms (cf Rickard \cite[Corollary 3.5]{Ri1}).

The result has many consequences. In particular
$$X\otimes_A^{\mathbb L}(-\otimes_A^{\mathbb L}Hom_A(X,A)):D^b(A\otimes_KA^{op})\lra
D^b(B\otimes_KB^{op})$$
is an equivalence. Therefore, if $D^b(A)\simeq D^b(B)$ there exists a two-sided
tilting complex realising an equivalence. The given equivalence one started with may be different. Then, this two-sided tilting complex induces the equivalence of the derived
categories of bimodules, and this then has the property that
$$End_{D^b(A\otimes A^{op})}(A)\stackrel{\simeq }{\lra} End_{D^b(B\otimes B^{op})}(B)$$
is an equivalence.

It is known that the for any algebra $C$ the module category is a full subcategory of the derived category
by identifying a module with the complex concentrated in a single degree $0$ (cf e.g.
Verdier \cite[Chapitre III Section 1.2.9]{Verdier}).
Therefore
for a $C$-module $M$ one gets $End_C(M)\simeq End_{D^b(C)}(M)$.
Hence $End_{D^b(A\otimes A^{op})}(A)\simeq End_{A\otimes A^{op}}(A)\simeq Z(A)$
and $End_{D^b(B\otimes B^{op})}(B)\simeq End_{B\otimes B^{op}}(B)\simeq Z(B)$.

\begin{Cor} (Rickard \cite[Proposition 2.5]{Ri1})
Let $K$ be a field and let $A$ and $B$ be finite dimensional $K$-algebras.
If $D^b(A)\simeq D^b(B)$, then a choice of a two-sided tilting complex $X$
realising this equivalence as tensor product induces an isomorphism $\varphi_X$
as algebras between the centre of $A$ and the centre of $B$.
\end{Cor}

We see that still an equivalence between the derived categories of finite dimensional
algebras yield an isomorphism between the centres of the algebras. The isomorphism
however is far less explicit and somewhat complicate.

Despite these difficulties we obtain the following result.

\begin{Theorem} \cite[Theorem 1]{Kuelsquest}\label{InvarianceKuelshammerideals}
Let $A$ and $B$ be finite dimensional symmetric $K$-algebras over a perfect
field $K$ of characteristic $p>0$. Suppose that $D^b(A)\simeq D^b(B)$
as triangulated categories. Then any choice of a two-sided tilting complex
$X$ yields an isomorphism $\varphi_X:Z(A)\lra Z(B)$ satisfying
$\varphi_X(T_n(A)^\perp)=T_n(B)^\perp$ for all $n\in\N$.
\end{Theorem}

It is worth writing that we use here the
mapping that is induced by a functor on the morphisms. We recall that possibly
non standard derived equivalences exist. Non standard derived equivalences
are not standard only on morphisms, but we only use morphisms here. Hence non standard
derived equivalences would possibly induce an isomorphism between the centres
which does not preserve the K\"ulshammer ideal structure.

The proof is much more involved than the proof for Morita invariance
in the sense that one needs to reformulate the
construction of K\"ulshammer ideals in a ''derived category
readable form''. Instead of
explicit constructions of particular sets one needs to argue via
homological properties
of morphism spaces.

One should mention that being symmetric is an invariant under derived equivalences.

\begin{Prop} (Rickard \cite[Corollary 5.3]{Ri1} for fields $R$, \cite{rogquest}
for more general rings)
Let $R$ be a Dedekind domain and let $A$ and $B$ be $R$-algebras
of finite rank over $R$ so that $D^b(A)\simeq D^b(B)$. Then if $A$ is symmetric,
$B$ is symmetric as well.
\end{Prop}

Applications of Theorem~\ref{InvarianceKuelshammerideals} will be given
in Section~\ref{applisection}. Actually the invariance of K\"uls\-hammer ideals
proved to be a rather powerful tool in particular in order to
distinguish algebras given by quivers and relations where the relations
depend on certain parameters. The structure of the quotient
$Z(A)/T_n(A)^\perp$ tends to depend on the
parameters in several cases.

\subsection{Stable invariance}

\label{stableinvariance}

As usual the stable category $A-\ul{mod}$ of a module category $A-mod$ is the category with
objects being $A$-modules and morphisms between two $A$-modules in the stable
category are equivalence classes of morphisms between these $A$-modules
modulo those which factor through a projective $A$-module. We denote
by $\ul{Hom}_A(M,N)$ the morphisms in $A-\ul{mod}$ from $M$ to $N$.

Equivalences between stable categories can behave badly in general.
An example was given by Auslander and Reiten in 1973 \cite[Example 3.5]{AR}.

\begin{Example}\label{ARindec}
Let $K$ be a field. Then the algebras
$$A:=\left(\begin{array}{cc}K&K\\0&K\end{array}\right)\times \left(\begin{array}{cc}K&K\\0&K\end{array}\right)$$
and
$$
B:=\left(\begin{array}{ccc}K&K&K\\0&K&K\\0&0&K\end{array}\right)\left/
\left(\begin{array}{ccc}0&0&K\\0&0&0\\0&0&0\end{array}\right)\right.
$$
the quotient of the upper triangular matrix ring by the ideal generated by the upper
right component matrices
have equivalent stable categories. Indeed,
$A$ has six indecomposable modules: four simple modules,
two of which are projective, and the projective cover of
the non projective simple.
Hence
the stable module category is equivalent to two copies of $K-mod$.
The algebra $B$ has five indecomposable modules: three simple modules,
one of which is projective and the projective covers of the other two
simple modules. Hence the
stable module category is equivalent to two copies of $K-mod$ as well.

The algebras $A$ and $B$ are hence stably equivalent and
$B$ is indecomposable whereas $A$ is not. Neither
$A$ nor $B$ has any simple direct factor.
\end{Example}

Given two self-injective algebras $A$ and $B$ and suppose $X$ is
a complex of $A-B$-bimodules inducing a standard equivalence
$D^b(B)\lra D^b(A)$. The quasi-inverse is again a standard equivalence,
given by a complex $Y$ of $B-A$-bimodules. Then a somewhat technical construction
on $X$ and on $Y$ produces an $A-B$-bimodule $M$, projective on either side, and a
$B-A$-bimodule $N$, projective on either side,
so that the $A-A$-bimodule $M\otimes_BN$ is isomorphic to $A\oplus P$ for some
projective $A-A$-bimodule $P$, and so that one has an isomorphism of $B-B$-bimodules
$N\otimes_AM\simeq B\oplus Q$
for some projective $B-B$-bimodule $Q$.
This motivated Brou\'e to define a class of stable equivalences with nicer properties.

\begin{Def} (Brou\'e \cite[Section 5]{Broue1994})
 Let $K$ be a commutative ring, let $A$ and $B$ be two $K$-algebras and
 let $M\in A\otimes_KB^{op}-mod$ and let $N\in B\otimes_KA^{op}-mod$.
 Then $(M,N)$ is said to induce a stable equivalence of Morita type if
 \begin{itemize}
 \item $M$ as well as $N$ are projective as $A$-modules and as $B$-modules.
 \item $M\otimes_BN\simeq A\oplus P$ as $A-A$-bimodules for a projective $A-A$-bimodule $P$
 \item $N\otimes_AM\simeq B\oplus Q$ as $B-B$-bimodules for a projective $B-B$-bimodule $Q$.
 \end{itemize}
\end{Def}

\begin{Rem}
Liu shows in \cite[Theorem 2.2]{Liu2008}
that a stable equivalence of Morita type between
two finite dimensional algebras with no separable summands restricts
to a stable equivalence between their summands. Therefore the algebras
$A$ and $B$ in Example~\ref{ARindec}
are not stably equivalent of Morita type.
\end{Rem}

In the meantime many properties have been shown to be invariant under
stable equivalence
of Morita type, whereas the general stable equivalences are still rather
poorly understood.

In particular, the following definition will be of importance in our discussion.
Let $A$ be a $K$-algebra. Then
$$Z(A)=Hom_{A\otimes_KA^{op}}(A,A)$$
and define the stable centre
$$Z^{st}(A)=\ul{Hom}_{A\otimes_KA^{op}}(A,A)$$
The natural homomorphism
$$Hom_{A\otimes_KA^{op}}(A,A)\lra \ul{Hom}_{A\otimes_KA^{op}}(A,A)$$
has a kernel denoted $Z^{pr}(A)$, the projective centre.

\begin{Prop} (Brou\'e \cite[Proposition 5.4]{Broue1994})
Let $A$ and $B$ be finite dimensional $K$-algebras and let $(M,N)$ be bimodules
inducing a stable equivalence of Morita type. Then
$Z^{st}(A)\simeq Z^{st}(B)$ as algebras.
\end{Prop}

How can we determine $Z^{pr}(A)$? This is a result of Liu, Zhou and the author.
The Cartan matrix of the algebra $A$ is denoted by $C_A$. Recall that the Cartan matrix
is square of size $n$, where $n$ is the number of simple $A$-modules up to
isomorphism. If one labels the rows and the columns by the isomorphism classes
$[S]$ of simple modules, then we have that the coefficient in position $([S],[T])$
is $Hom_A(P_S,P_T)$, where $P_S$ denotes the projective cover of $S$.
Hence $C_A$ has integer coefficients and can therefore
be interpreted as a linear endomorphism of the $K$-vector space $K^n$.

\begin{Prop} (Liu, Zhou, Zimmermann
\cite[Proposition 2.4, Corollary 2.9, Lemma 7.8, Proposition 7.10]{LZZ})
Let $K$ be an algebraically closed field and let $A$ be a
finite dimensional symmetric
$K$-algebra. Then $Z^{pr}(A)\subseteq Z(A)\cap soc(A)=R(A)$
and $dim_K(Z^{pr}(A))=rank_KC_A$.
\end{Prop}

In order to adapt the K\"ulshammer ideal theory for stable equivalences of Morita type
we need to replace the centre by the stable centre since
we know that a stable equivalence gives an isomorphism of the stable centres and
we do not have enough information about the centre.

Further we need to find a replacement of $A/[A,A]$. For this purpose we recall
the Hattori-Stallings trace which was generalised by Bouc to a trace function on the whole Hochschild homology.

\begin{Def} (Bouc~\cite{Bouc} for higher dimensional Hochschild homology,
Hattori-Stallings in degree $0$)\label{BoucDef}
Let $K$ be a field and let $A$ and $B$ be two finite dimensional $K$-algebras.
Given an $A-B$-bimodule $M$
which is projective as $B$-module, then there are elements
$m_i\in M,\varphi_i\in Hom_B(M,B)$, for $i=1,\dots,n$ so that
the identity on $M$ in $End_B(M)\simeq M\otimes_BHom_B(M,B)$ is mapped to
$\sum_{i=1}^nm_i\otimes\varphi_i$. The fact that $M$ is an $A-B$-bimodule
gives a mapping
\begin{eqnarray*}
A&\stackrel{\alpha_M}{\lra}&End_B(M)\simeq M\otimes_B Hom_B(M,B)\\
a&\mapsto&\sum_{i=1}^n(am_i)\otimes\varphi_i
\end{eqnarray*}
We produce $$eval:M\otimes_B Hom_B(M,B)\lra B/[B,B]$$
by $\mbox{eval}(m\otimes \psi):=\psi(m)+[B,B]$ for $\psi\in Hom_B(M,B)$ and $m\in M$.
The composition $\mbox{eval}\circ\alpha$ factorises through $A/[A,A]$ and the resulting
mapping $$A/[A,A]\lra B/[B,B]$$
is called the trace of $M$, denoted by $tr_M$. Similar statements hold if $M$ is
projective on the left.
\end{Def}

Using the Hattori-Stallings trace we define

\begin{Def}\cite[Defintion 4.1]{LZZ}
Let $A$ be a finite dimensional $K$-algebra
$$HH_0^{st}(A):=\bigcap_{P \mbox{ \scriptsize projective indecomposable } A-mod}ker(tr_P)$$
observing that any projective $A$-module $P$ is a $K-A$-bimodule as required
by Definition~\ref{BoucDef}.
\end{Def}

With this preparation we obtain that the dimension of $HH_0^{st}(A)$
is an invariant under stable equivalence of Morita type. Denote by $\ell(A)$
the number of simple $A$-modules up to isomorphism.

\begin{Theorem}\label{LZZtheorem}
Let $K$ be an algebraically closed field and let $A$ and $B$ be
finite dimensional
$K$-algebras without any semisimple direct factor
and suppose that $A$ and $B$ are stably equivalent of Morita type.
\begin{itemize}
\item (Liu, Zhou, Zimmermann \cite[Theorem 6.1]{LZZ})
Then
$$dim_K(A/[A,A])=dim_K(B/[B,B])\Leftrightarrow \ell(A)=\ell(B).$$
Moreover
$\mathrm{dim}(T_n(A)/K(A))=\mathrm{dim}(T_n(B)/K(B))$.
\item (Liu, Zhou, Zimmermann \cite[Corollary 6.2]{LZZ})
If in addition $A$ is symmetric then
$$\mathrm{dim}_K(Z(A))=\mathrm{dim}_K(Z(B))\Leftrightarrow \ell(A)=\ell(B).$$
and (K\"onig, Liu, Zhou \cite[Proposition 5.8]{KLZ})
$$Z(A)/T_n(A)^\perp\simeq Z(B)/T_n(B)^\perp$$
for all $n\geq 1$.
\end{itemize}
\end{Theorem}

The proof of the first part uses first that $HH_0^{st}(A)$ is an
invariant under stable equivalences of Morita type. Then one shows \cite[Theorem 4.4]{LZZ} that
$$dim(HH_0^{st}(A))+rank_K(C_A)=dim(A/[A,A])$$
Further it is shown in \cite[Section 5]{LZZ} that $rank_K(C_A)-dim_K(K\otimes_\Z K_0(A))$ equals the dimension
of the so-called stable Grothendieck group, which is known to be an invariant under stable
equivalences of Morita type by work of Xi~\cite[Section 5]{Xi2008}.

We should mention the long standing Auslander-Reiten conjecture.

\begin{Conj} (Auslander-Reiten \cite[page 409 Conjecture 5]{ARS}) \label{ARConj}
Let $A$ and $B$ be finite dimensional $K$-algebras. If $A$ and $B$ are stably equivalent,
then the number of simple non projective $A$-module up to isomorphism equals the number of
simple non projective $B$-modules up to isomorphism.
\end{Conj}

A priori I feel that there is no obvious reason why the invariance of
the number of simple $A$-modules has anything to do with the
commutator quotient. This fact appears somewhat surprisingly in this connection.

The conjecture has been verified for quite a few classes of algebras. Few general
positive results are known so far.

\section{Applications}
\label{applisection}

In the last five years
K\"ulshammer ideals were successfully employed to distinguish algebras up to
derived and up to stable equivalences of Morita type for various classes
of algebras which were extremely difficult to deal with previously. In particular if
two algebras are defined by the same quiver $Q$ and a set of relations $I(c)$
subject to some parameter $c$ in the base field, then the technique of
computing K\"ulshammer ideals and the quotient of the centre by the ideal proves to be
fruitful in various cases.

\subsection{Algebras of dihedral, semidihedral and quaternion type}

Many authors during the last decades
proved ring theoretic properties for group algebras, and still the
question is an active field of research. In particular
many properties are shown to hold for the Cartan matrices
and the occurrence of certain components in the
stable Auslander-Reiten quiver for blocks of group algebras
with dihedral, semidihedral or quaternion defect groups. Moreover
it was shown at that time that a block of a group algebra is of tame representation
type (cf Section~\ref{polygrowthsection} below for the precise definition)
if and only if the defect group is a dihedral, a semidihedral or a quaternion
group.

Erdmann showed in \cite{ErdmannLNM} that these properties
determine the Morita equivalence classes of these algebras
as belonging to a finite number of families, given by quivers with relations,
subject to certain parameters in the relations.
Up to these parameters in the relations the algebras are classified in
a finite number of classes up to Morita equivalences.

Holm  \cite{Holmhabil} classified further these Morita equivalence classes
up to derived equivalences. Many classes merge to a common derived equivalence class.
However, Holm could not determine for
a certain number of parameters if two algebras within
one class but with different
parameters are derived equivalent or not.

K\"ulshammer ideals manage to distinguish derived equivalence classes in some of these cases.

We display Thorsten Holm's list \cite{Holmhabil} of
algebras of dihedral, semidihedral and quaternion type up to
derived equivalences. Each of these types form a family.
Each family is subdivided into three
subclasses: algebras with one simple module, algebras with two
simple modules and algebras with three simple modules. Each
subfamily contains
algebras defined by quivers and relations, depending on parameters.

{\small $$\begin{array}{c||c|c|c}
&\mbox{dihedral}&\mbox{semidihedral}&\mbox{quaternion}\\ \hline
\mbox{1 simple}&K[X,Y]/(XY,X^m-Y^n),&SD(1\mathcal{A})_1^k,
k\geq 2;&Q(1\mathcal{A})_1^k, k\geq 2;\\
& m\geq n\geq 2, m+n>4;&&\\
&&&\\
&D(1\mathcal{A})_1^{1}=K[X,Y]/(X^2,Y^2); &&  \\
&&&\\
& (char K=2)   &(char K=2)\ SD(1\mathcal{A})_2^k(c,d)
&(char K=2)\ Q(1\mathcal{A})_2^k(c,d),\\
&K[X,Y]/(X^2,YX-Y^2);&k\geq 2, (c,d)\neq (0,0);&k\geq 2, (c,d)\neq (0,0);\\
&&&\\
&D(1\mathcal{A})^k_1, k\geq 2;&  & \\
&&&\\
&&&\\
&(char K=2)\ D(1\mathcal{A})^k_2(d),&&\\
&k\geq 2, d=0 \ or\  1;& &\\ \hline
\mbox{2 simples}&D(2\mathcal{B})^{k,s}(c),
&SD(2\mathcal{B})^{k,t}_1(c)&Q(2\mathcal{B})_1^{k,s}(a,c)\\
&k\geq s\geq 1, c\in\{0,\  1\}&k\geq 1, t\geq 2, c\in\{0, \ 1\};
&k\geq 1, s\geq 3, a\neq 0; \\
&&&\\
&&SD(2\mathcal{B})^{k,t}_2(c)&\\
& &k\geq 1, t\geq 2,  &\\
& &  k+t\geq 4,  c\in\{0, \ 1\};&\\
\hline
\mbox{3 simples}&D(3\mathcal{K})^{a,b,c},
&SD(3\mathcal{K})^{a,b,c}&Q(3\mathcal{K})^{a,b,c}\\
&a\geq b\geq c\geq 1;& a\geq b\geq c\geq 1, a\geq 2;
&a\geq b\geq c\geq 1, b\geq 2,\\
&&& (a, b, c)\neq (2,2,1);\\
&&&\\
 &D(3\mathcal{R})^{k,s,t,u},&&Q(3\mathcal{A})^{2,2}_1(d)\\
 &s\geq t\geq u\geq k\geq 1, t\geq 2&&d\not\in\{ 0, 1\}\\
\hline
\end{array}$$}

All algebras with one simple module in the above list have the
quiver of type $1\mathcal{A}$
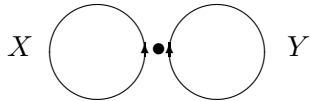
$$\unitlength0.6cm
\begin{picture}(7,3)
\put(5.45,2){\circle{2.0}}
\put(2.8,2){\circle{2.0}}
\put(3.85,1.9){\vector(0,1){0.3}}
\put(4,1.9){$\bullet$}
\put(4.4,1.9){\vector(0,1){0.3}}
\put(0.8,1.9){$X$} \put(7,1.9){$Y$}
\end{picture}$$

The quivers of the algebras of type $2\mathcal{B}$, $3\mathcal{K}$,
$3\mathcal{A}$ and $3\mathcal{R}$ are respectively:

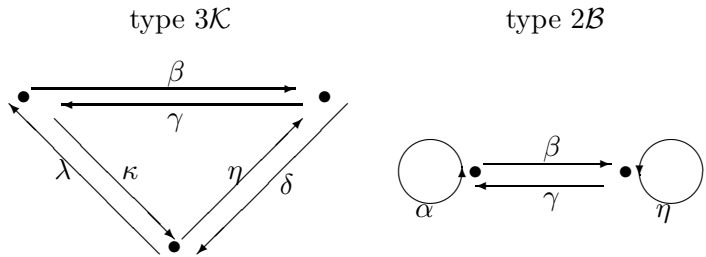
\unitlength1cm
\begin{picture}(14,5)
\put(3.5,4){type $3\mathcal{K}$} \put(2,3){$\bullet$}
\put(4,1){$\bullet$} \put(6,3){$\bullet$}
\put(2.2,3.2){\vector(1,0){3.5}}
\put(5.8,2.99){\vector(-1,0){3.2}}
\put(2.5,2.8){\vector(1,-1){1.6}} \put(3.9,1){\vector(-1,1){2}}
\put(4.2,1.2){\vector(1,1){1.6}} \put(6.4,3){\vector(-1,-1){2}}
\put(4,3.3){$\beta$} \put(4,2.7){$\gamma$} \put(5.5,1.8){$\delta$}
\put(4.8,2){$\eta$} \put(2.5,2){$\lambda$} \put(3.4,2){$\kappa$}

\put(8.5,4){type $2\mathcal{B}$} \put(8,2){$\bullet$}
\put(10,2){$\bullet$} \put(8.2,2.2){\vector(1,0){1.7}}
\put(9.8,1.9){\vector(-1,0){1.7}} \put(7.5,2.1){\circle{.8}}
\put(10.7,2.1){\circle{.8}} \put(7.92,2.15){\vector(0,1){.01}}
\put(10.28,2.05){\vector(0,-1){.01}} \put(7.3,1.5){$\alpha$}
\put(10.5,1.5){$\eta$} \put(9,2.3){$\beta$} \put(9,1.65){$\gamma$}

\end{picture}

\bigskip

\unitlength1cm
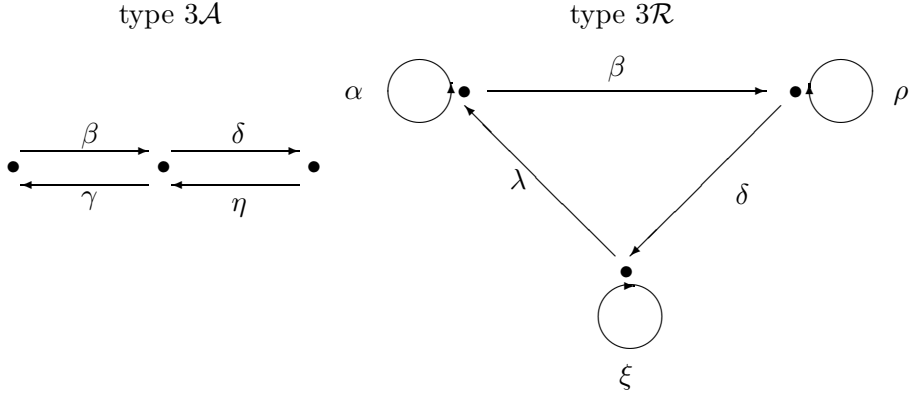
\noindent
\begin{picture}(15,6)

\put(3.5,6){type $3\mathcal{A}$} \put(2,4){$\bullet$}
\put(4,4){$\bullet$} \put(6,4){$\bullet$}
\put(2.2,4.3){\vector(1,0){1.7}}
\put(3.9,3.85){\vector(-1,0){1.7}}
\put(4.2,4.3){\vector(1,0){1.7}}
\put(5.9,3.85){\vector(-1,0){1.7}}
\put(5,3.5){$\eta$} \put(3,4.4){$\beta$} \put(3,3.6){$\gamma$}
\put(5,4.4){$\delta$}

\put(9.5,6){type $3\mathcal{R}$}
\put(8,5){$\bullet$}
\put(10.15,2.6){$\bullet$}
\put(12.4,5){$\bullet$}
\put(12.68,5.15){\vector(0,1){.05}}
\put(8.4,5.1){\vector(1,0){3.7}}
 \put(10.1,2.9){\vector(-1,1){2}}
 \put(12.3,4.9){\vector(-1,-1){2}}
  \put(10,5.3){$\beta$}
\put(11.7,3.6){$\delta$}
\put(8.7,3.8){$\lambda$}

\put(7.5,5.1){\circle{.8}}\put(7.92,5.2){\vector(0,1){.01}}
\put(13.1,5.1){\circle{.8}}\put(12.7,53.2){\vector(0,1){.01}}
\put(10.3,2.1){\circle{.8}}\put(10.35,2.5){\vector(1,0){.01}}
\put(6.5,5){$\alpha$} \put(10.15,1.2){$\xi$} \put(13.8,5){$\rho$}

\end{picture}

The relations are respectively
\begin{eqnarray*}
D(1\mathcal{A})^k_1 &: &  X^2, Y^2,(XY)^k-(YX)^k;\\
D(1\mathcal{A})^k_2(d)&: &  X^2-(XY)^k,Y^2-d\cdot(XY)^k,(XY)^k-(YX)^k,(XY)^kX,(YX)^kY;\\
SD(1\mathcal{A})_1^k&: &(XY)^k-(YX)^k,(XY)^kX,Y^2,X^2-(YX)^{k-1}Y;\\
SD(1\mathcal{A})_2^k(c,d)&:&(XY)^k-(YX)^k,(XY)^kX,Y^2-d(XY)^k,\\
&&\phantom{K<X,Y>/}X^2-(YX)^{k-1}Y+c(XY)^k;\\
Q(1\mathcal{A})_1^k&:&(XY)^k-(YX)^k,(XY)^kX,Y^2-(XY)^{k-1}X,X^2-(YX)^{k-1}Y;\\
Q(1\mathcal{A})_2^k(c,d)&:&X^2-(YX)^{k-1}Y-c(XY)^k,Y^2-(XY)^{k-1}X-d(XY)^k,\\
&&\phantom{K<X,Y>/}(XY)^k-(YX)^k,(XY)^kX,(YX)^kY.
\end{eqnarray*}
as well as
\begin{eqnarray*}
D(2\mathcal{B})^{k, s}(c) &: &
\beta\eta, \eta\gamma, \gamma\beta, \alpha^2-c(\alpha\beta\gamma)^k,
(\alpha\beta\gamma)^k-(\beta\gamma\alpha)^k, \eta^s-(\gamma\alpha\beta)^k;\\
SD(2\mathcal{B})^{k, t}_1(c)&: & \gamma\beta, \eta\gamma,
\beta\eta,
\alpha^2-(\beta\gamma\alpha)^{k-1}\beta\gamma-c(\alpha\beta\gamma)^k,
\eta^t-(\gamma\alpha\beta)^k,
(\alpha\beta\gamma)^k-(\beta\gamma\alpha)^k;\\
SD(2\mathcal{B})^{k, t}_2(c)&:&
\beta\eta-(\alpha\beta\gamma)^{k-1}\alpha\beta,
\eta\gamma-(\gamma\alpha\beta)^{k-1}\gamma\alpha,
\gamma\beta-\eta^{t-1}, \alpha^2-c(\alpha\beta\gamma)^k,
\beta\eta^2, \eta^2\gamma;\\
Q(2\mathcal{B})^{k, s}_1(a,c)&:&\gamma\beta-\eta^{s-1},
\beta\eta-(\alpha\beta\gamma)^{k-1}\alpha\beta,
\eta\gamma-(\gamma\alpha\beta)^{k-1}\gamma\alpha,
 \\
&&\phantom{K<X,Y>/}\alpha^2-a(\beta\gamma\alpha)^{k-1}\beta\gamma-c
(\beta\gamma\alpha)^k, \alpha^2\beta, \gamma\alpha^2;
\end{eqnarray*}
\begin{eqnarray*}
D(3\mathcal{K})^{a, b, c} &: & \beta\delta, \delta\lambda,
\lambda\beta, \gamma\kappa, \kappa\eta, \eta\gamma,
(\beta\gamma)^a-(\kappa\lambda)^b,
(\lambda\kappa)^b-(\eta\delta)^c, (\delta\eta)^c-(\gamma\beta)^a;\\
D(3\mathcal{R})^{k, s, t, u} &: & \alpha\beta, \beta\rho,
\rho\delta, \delta\xi, \xi\lambda, \lambda\alpha,
\alpha^s-(\beta\delta\lambda)^k, \rho^t-(\delta\gamma\beta)^k,
\xi^u-(\lambda\beta\delta)^k  ;\\
SD(3\mathcal{K})^{a, b, c} &: & \kappa\eta, \eta\gamma, \gamma\kappa,
\delta\gamma-(\gamma\alpha)^{a-1}\gamma,
\beta\delta-(\kappa\lambda)^{b-1}\kappa, \lambda\beta-(\eta\delta)^{c-1}\eta;\\
Q(3\mathcal{K})^{a, b, c} &: &
\beta\delta-(\kappa\lambda)^{a-1}\kappa,
\eta\gamma-(\lambda\kappa)^{a-1}\lambda,
\delta\lambda-(\gamma\beta)^{b-1}\gamma,
\kappa\eta-(\beta\gamma)^{b-1}\beta,\\
&&\phantom{K<X,Y>/}
\lambda\beta-(\eta\delta)^{c-1}\eta,  \gamma\kappa-(\delta\eta)^{c-1}\delta,
\gamma\beta\delta, \delta\eta\gamma, \lambda\kappa\eta;\\
Q(3\mathcal{A})^{2,2}_1(d) &: & \beta\delta\eta-\beta\gamma\beta,
\delta\eta\gamma-\gamma\beta\gamma,
\eta\gamma\beta-d\eta\delta\eta,
\gamma\beta\delta-d\delta\eta\delta, \beta\delta\eta\delta,
\eta\gamma\beta\gamma.
\end{eqnarray*}

For the dihedral type algebras with two simple modules a result of
Kauer and Roggenkamp \cite[Corollary 5.3]{KR}
show that the parameters $c=0$ and $c=1$ yield different derived equivalence
classes of algebras. The method employed there is rather involved. The authors define
graph algebras and show that being a graph algebra is invariant under
derived equivalences. Further for one of the scalars the algebra is a
graph algebra, for the other it is not. Holm and the author
gave a much simpler proof in  \cite{hz-tame} avoiding graph algebras.

The semidihedral type case can be dealt with at least partially.
Again we consider the derived equivalence classes of semidihedral
type algebras  with two simple modules, and again the question
if the parameters $c=0$ and $c=1$ yield different derived equivalence
classes was open.

\begin{Theorem}(Holm and Zimmermann
\cite[Theorem 1.1, Theorem 1.2 and Theorem 1.3]{hz-tame}) \label{hz-tame2}
\label{thm1-intro-semidihedral}
Let $K$ be an algebraically closed field of characteristic 2.
\begin{itemize}
\item
For any given integers $k,s\ge 1$
consider the algebras of dihedral type
$D(2A)^{k,s}(c)$ for the scalars $c=0$ and $c=1$.
Suppose that $k\geq 2$.
Suppose if $k=2$ then $s\geq 3$ is odd, and if $s=2$ then
$k\geq 3$ is odd.

Put $A_0^{k,s}:=D(2A)^{k,s}(0)$ and $A_1^{k,s}:=D(2A)^{k,s}(1)$.

Then  the factor rings $Z(A^{k,s}_0)/T_1(A^{k,s}_0)^{\perp}$
and $Z(A^{k,s}_1)/T_1(A^{k,s}_1)^{\perp}$ are not isomorphic as
rings.

In particular, the algebras $D(2A)^{k,s}(0)$ and $D(2A)^{k,s}(1)$
are not derived equivalent and are not stably equivalent of Morita type.

\item
For any given integers $k\geq 1$ and $s\geq 1$, consider the algebras of
semidihedral type
$SD(2B)_1^{k,s}(c)$ for the scalars $c=0$ and $c=1$.
Suppose that $k\geq 2$.
Suppose that if $k=2$ then $s\geq 3$ is odd, and if $s=2$ then
$k\geq 3$ is odd.

Put $B_0^{k,s}:=SD(2B)_1^{k,s}(0)$ and $B_1^{k,s}:=SD(2B)_1^{k,s}(1)$.

Then
the factor rings $Z(B^{k,s}_0)/T_1(B^{k,s}_0)^{\perp}$
and $Z(B^{k,s}_1)/T_1(B^{k,s}_1)^{\perp}$ are not isomorphic as
rings.

In particular, the algebras $SD(2B)_1^{k,s}(0)$ and $SD(2B)_1^{k,s}(1)$
are not derived equivalent and are not stably equivalent of
Morita type in these cases.

\item
For any given integers $k\geq 1$ and $s\geq 1$, consider the algebras of
semidihedral type
$SD(2B)_2^{k,s}(c)$ for the scalars $c=0$ and $c=1$.
Suppose $k\geq 2$.
Put $C_0^{k,s}:=SD(2B)_2^{k,s}(0)$ and $C_1^{k,s}:=SD(2B)_2^{k,s}(1)$.

If the parameters $k$ and $s$ are both odd, then
the factor rings $Z(C^{k,s}_0)/T_1(C^{k,s}_0)^{\perp}$
and $Z(C^{k,s}_1)/T_1(C^{k,s}_1)^{\perp}$ are not isomorphic as
rings.

In particular, the algebras $SD(2B)_2^{k,s}(0)$ and $SD(2B)_2^{k,s}(1)$
are not derived equivalent and are not
stably equivalent of Morita type in these cases.
\end{itemize}
\end{Theorem}

\begin{Rem}
\begin{enumerate}
\item
We should mention that actually the dimension of the quotients
of the centres modulo the K\"ulshammer ideals do not depend on the
scalar $c$. The algebra structure of the quotient is needed.

\item
It is worth noticing that these algebras are all symmetric and so \cite{Kuelsquest}
applies directly. Moreover, the dimension of the centre of the algebra
equals the dimension of the
quotient of the algebra by the commutator subspace. This immediate
consequence of the
fact that the algebras are symmetric is not clear in case the algebra
is selfinjective only.
Example~\ref{twistedcentreexample} gives an easy example for how
complicated the situation might become
already for very small selfinjective algebras. For general selfinjective
algebras a rather sophisticated theory needs to be developed in order to
compute the commutator subspace.
\end{enumerate}
\end{Rem}

\begin{Rem}
I would like to mention that in
\cite[Theorem 1.1]{hz-tame}
the condition that $k\geq 2$
in case of algebras of dihedral type in Theorem~\ref{hz-tame2} is unfortunately missing.
The condition is necessary. A recent result
of Frauke Bleher \cite[Theorem 2]{Bleher} determines the parameter $c$ in the relations
for a specific group by completely different methods. If $k=1$ would be
allowed, then the parameter would be different than determined by Bleher.
This observation is due to Zhou.
\end{Rem}

In recent work Zhou and the author studied if the derived equivalence
classification of Holm
of dihedral, semidihedral and quaternion type algebras also gives a
classification
up to stable equivalence. Some partial statements are already given in Theorem~\ref{hz-tame2}.

One has to deal with several additional problems for stable
equivalences of Morita type.

The first problem is that
the Auslander-Reiten conjecture is open, that is we might a priori have
a stable equivalence of Morita type between two algebras of dihedral,
semidihedral or quaternion
type with different number of simple modules. This does not happen for
derived equivalences since there the rank of the Grothendieck group is an invariant.

The second problem is that
a derived equivalence between an algebra and a local algebra is in fact a Morita equivalence. This
was shown by Roggenkamp and the author~\cite[Section 5]{Mexico}. The statement is
false for stable equivalences of Morita type. Given a finite group $G$ and
a field $K$ of characteristic $p$ dividing the order of $G$
a $KG$-module $M$ is endotrivial if the $KG$-module $End_K(M)$
has the property $End_K(M)\simeq K\oplus P$ for $K$ being the
trivial $K$-module and $P$ a projective $KG$-module.
Every endotrivial module over a $p$-group gives a stable
self-equivalence of Morita type for the group ring over this $p$-group. The set
of endotrivial modules over a fixed $p$-group up to some equivalence relation
form a group, whose structure was completely determined by Carlson and Th\'evenaz,
and which is non trivial free abelian in most cases \cite{CarlsonThevenaz}.

\medskip

Since the statement of the result might be technical
for the non specialist reader we illustrate the result in a coarser form.

Recall that we have the following rough classification of algebras up to derived equivalences.

{\small $$\begin{array}{c||c|c|c}
&\mbox{dihedral}&\mbox{semidihedral}&\mbox{quaternion}\\ \hline
\mbox{1 simple}&\mbox{ five types of algebras }&\mbox{ two types of algebras }&\mbox{ two types of algebras }\\
& \mbox{ depending on parameters }&\mbox{ depending on parameters }&\mbox{ depending on parameters }\\
\hline
\mbox{2 simples}&\mbox{ one type of algebras }
&\mbox{ two types of algebras }&\mbox{ one type of algebras }\\
&\mbox{ depending on parameters }&\mbox{ depending on parameters }
&\mbox{ depending on parameters } \\
\hline
\mbox{3 simples}&\mbox{ two types of algebras }
&\mbox{ one type of algebras }&\mbox{ two types of algebras }\\
&\mbox{ depending on parameters }&\mbox{ depending on parameters }
&\mbox{ depending on parameters }\\
\hline
\end{array}$$}

Theorem~\ref{main} below states mainly that the columns are
preserved under stable equivalences
of Morita type and that the rows are preserved under
stable equivalences of Morita type.

The actual statement is finer than this, but this scheme
gives a relatively good approximation
of what is proved in Theorem~\ref{main}.

The details we obtain are given in the following result.

\begin{Theorem} (Zhou and Zimmermann \cite[Theorem 7.1]{Modulestructure})
\label{main}
Let $K$ be an algebraically closed field.

Suppose $A$ and $B$ are indecomposable algebras which
are stably equivalent of Morita type.

\begin{itemize}
\item
If $A$ is an algebra of dihedral type, then $B$ is of
dihedral type. If $A$ is of semidihedral type, then
$B$ is of semidihedral type. If $A$ is of quaternion type then
$B$ is of quaternion type.

\item If $A$ and $B$ are of dihedral, semidihedral or quaternion type, then
$A$ and $B$ have the same number of simple modules.

\item Let $A$ be an algebra of
dihedral type.

\begin{enumerate}
\item
If $A$ is local, then $A$ is stably
equivalent of Morita type to one and exactly one  algebra in the
following list:
\begin{itemize}
\item $A_1(n, m)$ with $m\geq n\geq 2$ and $m+n>4$;
\item $C_1 $;
\item $D(1\mathcal{A})_1^k $ with $k\geq 2$;
\item  if $p=2$, $B_1$  and $D(1\mathcal{A})_2^k(d)$
with $k\geq 2$ and $d\in \{0, 1\}$, except that we don't
  know whether $D(1\mathcal{A})^k(0) $ and  $D(1\mathcal{A})^k(1)$ are
  stably equivalent of Morita type or not.
\end{itemize}

\item
If $A$ has two simple modules,
then $A$ is stably  equivalent of Morita type  to one and exactly
one of the following algebras:
$D(2\mathcal{B})^{k, s}(0)$ with $k\geq s\geq 1$ or if $p=2$, $D(2\mathcal{B})^{k, s}(1)$
with $k\geq s\geq 1$.

\item
If $A$ has three simple modules then
$A$ is stably  equivalent of Morita type  to one and exactly one of
the following algebras:  $D(3\mathcal{K})^{a, b, c}$ with $a\geq b\geq c\geq 1$
or   $D(3\mathcal{R})^{k, s, t, u}$ with $  s\geq t\geq u\geq k\geq 1$ and $t\geq 2$.
\end{enumerate}

\item Let $A$ be an algebra of semidihedral type.
\begin{enumerate}
\item If $A$ has one simple module then $A$ is stably equivalent
of Morita type to one of the following algebras:
$SD(1\mathcal{A})_1^k$ for $k\geq 2$ or
$SD(1\mathcal{A})_2^k(c,d)$ for $k\geq 2$ and $(c,d)\neq (0,0)$ if
the characteristic of $K$ is $2$. Different parameters $k$ yield
algebras in different stable
 equivalence classes of Morita type.

\item
If $A$ has two simple modules then $A$ is stably equivalent of Morita type to
$SD(2\mathcal B)_1^{k,s}(c)$ for $k\geq 1,s\geq 2,c\in\{0,1\}$
or to
$SD(2\mathcal B)_2^{k,s}(c)$ for $k\geq 1,s\geq 2,c\in\{0,1\},k+s\geq 4.$

\item
If $A$ has three simple modules, then $A$ is stably equivalent of Morita type to
one and only one algebra of the type $SD(3\mathcal{K})^{a,b,c}$ for
$a\geq b\geq c\geq 1$.
\end{enumerate}

\item Let $A$ be an algebra of quaternion type.

\begin{enumerate}
\item If $A$ has one simple modules, then $A$ is stably equivalent
of Morita type to one of the algebras $Q(1{\mathcal A})_1^k$ for
$k\geq 2$ or $Q(1{\mathcal A})_2^k(c,d)$ for $k\geq 2, (c,d)\neq
(0,0)$ if characteristic if the $K$ is $2$. Different parameters
$k$ yield algebras in different
stable equivalence classes of Morita type.

\item If $A$ has two simple modules then $A$ is stably equivalent of Morita
type to one of the algebras $Q(2{\mathcal B})_1^{k,s}(a,c)$
for $k\geq 1,s\geq 3,a\neq 0$.

\item If $A$ has three simple modules, then $A$ is stably equivalent of Morita type to
one of the algebras $Q(3{\mathcal K})^{a,b,c}$ for
$a\geq b\geq c\geq 1, b\geq 2, (a,b,c)\neq (2,2,1)$ or
$Q(3{\mathcal A})_1^{2,2}(d)$ for $d\in K\setminus\{0,1\}$.
Different parameters $a,b,c$ yield algebras in different stable equivalence
classes of Morita type.
\end{enumerate}

\end{itemize}
\end{Theorem}

One particularly nice consequence should be mentioned though.

\begin{Cor} (Zhou and Zimmermann \cite[Corollary 7.3]{Modulestructure})
The Auslander-Reiten conjecture \ref{ARConj} is true for algebras of
dihedral, semidihedral or quaternion type.
\end{Cor}

\subsection{Bia\l kowski-Erdmann-Skowro\'nski deformation of preprojective algebras}

\label{BESdeformedpreprojective}

Recently
Bia\l kowski, Erdmann and Skowro\'nski classified in \cite{BES} all selfinjective algebras
with the property that for all simple modules $S$ the third syzygy of
$S$ is again isomorphic
to $S$.  A recent survey on the circle around these questions was given by
Erdmann and Skowro\'nski in  \cite{periodic}.

The problem of classifying algebras so that $\Omega^2(S)\simeq S$ for all simple
modules was completely solved before
and the next most interesting case is
$\Omega^3(S)\simeq S$ for all simple
modules. In order to formulate the result \cite{BES} of Bia\l kowski, Erdmann and
Skowro\'nski we need to introduce deformed preprojective algebras as defined in \cite{BES}.

The preprojective algebra of type $A_n$ is given by the quiver

\unitlength 1cm
\begin{center}
\begin{picture}(10,2)
\put(1,1){$\bullet$}
\put(1.1,1.3){\vector(1,0){.8}}
\put(1.9,.9){\vector(-1,0){.8}}
\put(1.5,1.5){$a_1$}
\put(1.5,.5){$\ol a_1$}
\put(2,1){$\bullet$}
\put(2.1,1.3){\vector(1,0){.8}}
\put(2.9,.9){\vector(-1,0){.8}}
\put(2.5,1.5){$a_2$}
\put(2.5,.5){$\ol a_2$}
\put(3,1){$\bullet$}
\put(3.1,1.3){\vector(1,0){.8}}
\put(3.9,.9){\vector(-1,0){.8}}
\put(3.5,1.5){$a_3$}
\put(3.5,.5){$\ol a_3$}
\put(4,1){$\bullet$}
\put(4.1,1.3){\vector(1,0){.8}}
\put(4.9,.9){\vector(-1,0){.8}}
\put(4.5,1.5){$a_4$}
\put(4.5,.5){$\ol a_4$}
\put(5.5,1){$\dots\dots$}

\put(7.1,1.3){\vector(1,0){.8}}
\put(7.9,.9){\vector(-1,0){.8}}
\put(7.3,1.5){$a_{n-1}$}
\put(7.3,.5){$\ol a_{n-1}$}
\put(8,1){$\bullet$}
\end{picture}
\end{center}

subject to the relations
$$a_1\ol a_1=\ol a_{n-1}a_{n-1}=0\mbox{ and }
\ol a_ia_i=a_{i+1}\ol a_{i+1}\forall\;i\in\{1,2,\dots,n-2\}.
$$

The deformed preprojective algebra of type $D_{n+1}$ is given by the quiver

\unitlength 1cm
\begin{center}
\begin{picture}(10,3)
\put(1,0){$\bullet$}
\put(1.1,.3){\vector(1,1){.8}}
\put(1.9,.9){\vector(-1,-1){.8}}
\put(1,.7){$a_1$}
\put(1.5,.2){$\ol a_1$}
\put(1,2){$\bullet$}
\put(1.1,1.9){\vector(1,-1){.8}}
\put(1.9,1.3){\vector(-1,1){.8}}
\put(1.5,1.8){$\ol a_0$}
\put(.99,1.5){$a_0$}
\put(2,1){$\bullet$}
\put(2.1,1.3){\vector(1,0){.8}}
\put(2.9,.9){\vector(-1,0){.8}}
\put(2.5,1.5){$a_2$}
\put(2.5,.5){$\ol a_2$}
\put(3,1){$\bullet$}
\put(3.1,1.3){\vector(1,0){.8}}
\put(3.9,.9){\vector(-1,0){.8}}
\put(3.5,1.5){$a_3$}
\put(3.5,.5){$\ol a_3$}
\put(4,1){$\bullet$}
\put(4.1,1.3){\vector(1,0){.8}}
\put(4.9,.9){\vector(-1,0){.8}}
\put(4.5,1.5){$a_4$}
\put(4.5,.5){$\ol a_4$}
\put(5.5,1){$\dots\dots$}

\put(7.1,1.3){\vector(1,0){.8}}
\put(7.9,.9){\vector(-1,0){.8}}
\put(7.3,1.5){$a_{n-1}$}
\put(7.3,.5){$\ol a_{n-1}$}
\put(8,1){$\bullet$}
\end{picture}
\end{center}

subject to the relations
$$a_0\ol a_0=a_1\ol a_1=\ol a_{n-1}a_{n-1}=\ol a_1a_1+\ol a_0a_0+a_2\ol a_2+
f(\ol a_0a_0,\ol a_1a_1)=(\ol a_1a_1+\ol a_0a_0)^{n-2}=0$$
$$\mbox{ and }
\ol a_ia_i=a_{i+1}\ol a_{i+1}\forall\;i\in\{2,\dots,n-2\}.
$$
for some element $$f(X,Y)\in rad^2(K<X,Y>/(X^2,Y^2,(X+Y)^{n-1}).$$

The algebra of type $L_n$ for $n\geq 2$ was already displayed in
Example~\ref{ExampleFormIsNotSymmetric} and is given by the quiver
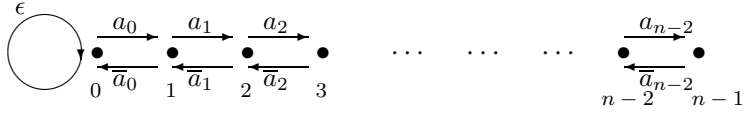
\unitlength1cm
\begin{center}
\begin{picture}(10,2)
\put(1,1){$\bullet$}
\put(1,0.5){\scriptsize $0$}
\put(2,1){$\bullet$}
\put(2,0.5){\scriptsize $1$}
\put(3,1){$\bullet$}
\put(3,0.5){\scriptsize $2$}
\put(4,1){$\bullet$}
\put(4,0.5){\scriptsize $3$}
\put(8,1){$\bullet$}
\put(7.8,0.4){\scriptsize $n-2$}
\put(9,1){$\bullet$}
\put(9,0.4){\scriptsize $n-1$}
\put(5,1){$\cdots$}
\put(6,1){$\cdots$}
\put(7,1){$\cdots$}
\put(1.1,1.3){\vector(1,0){.8}}
\put(2.1,1.3){\vector(1,0){.8}}
\put(3.1,1.3){\vector(1,0){.8}}
\put(8.1,1.3){\vector(1,0){.8}}
\put(8.9,.9){\vector(-1,0){.8}}
\put(3.9,.9){\vector(-1,0){.8}}
\put(2.9,.9){\vector(-1,0){.8}}
\put(1.9,.9){\vector(-1,0){.8}}
\put(0.4,1.1){\circle{1}}
\put(.88,1.1){\vector(0,-1){.1}}
\put(0,1.6){\small $\epsilon$}
\put(1.3,1.4){\small $a_0$}
\put(2.3,1.4){\small $a_1$}
\put(3.3,1.4){\small $a_2$}
\put(8.3,1.4){\small $a_{n-2}$}
\put(1.3,.65){\small $\ol a_0$}
\put(2.3,.65){\small $\ol a_1$}
\put(3.3,.65){\small $\ol a_2$}
\put(8.3,.65){\small $\ol a_{n-2}$}
\end{picture}
\end{center}
subject to the relations
$$a_i\overline a_i+\overline a_{i-1}a_{i-1}=0\mbox{ for all }
i\in\{1,\dots ,n-2\}~,
$$
$$\overline a_{n-2}a_{n-2}=0~,~~~~~
\epsilon^{2n}=0~,~~~~~
\epsilon^2+a_0\overline a_0+\epsilon^{3}p(\epsilon)=0
$$
for a polynomial $p(X)\in K[X]$. Denote by $L_n^p$ the deformed preprojective algebra
of type $L$ with deformation polynomial $p(X)$ and abbreviate $L_n^j:=L_n^{X^{2j}}$
for simplicity when no confusion may occur.

The deformed preprojective algebra of type $E_n$ for $n\in\{6,7,8\}$ is given by the
quiver
\unitlength 1.5cm
\begin{center}
\begin{picture}(10,3)
\put(1,1){$\bullet$}
\put(1.1,1.3){\vector(1,0){.8}}
\put(1.9,.9){\vector(-1,0){.8}}
\put(1.5,1.5){$a_0$}
\put(1.5,.5){$\ol a_0$}
\put(2,1){$\bullet$}
\put(2.1,1.3){\vector(1,0){.8}}
\put(2.9,.9){\vector(-1,0){.8}}
\put(2.3,1.5){$a_1$}
\put(2.5,.5){$\ol a_1$}
\put(3,1){$\bullet$}
\put(2.89,1.4){\vector(0,1){.8}}
\put(3.1,2.2){\vector(0,-1){.8}}
\put(2.6,1.8){$a_2$}
\put(3.2,1.8){$\ol a_2$}
\put(2.95,2.2){$\bullet$}

\put(3.1,1.3){\vector(1,0){.8}}
\put(3.9,.9){\vector(-1,0){.8}}
\put(3.5,1.5){$a_3$}
\put(3.5,.5){$\ol a_3$}
\put(4,1){$\bullet$}
\put(4.1,1.3){\vector(1,0){.8}}
\put(4.9,.9){\vector(-1,0){.8}}
\put(4.5,1.5){$a_4$}
\put(4.5,.5){$\ol a_4$}
\put(5.5,1){$\dots\dots$}

\put(7.1,1.3){\vector(1,0){.8}}
\put(7.9,.9){\vector(-1,0){.8}}
\put(7.3,1.5){$a_{n-2}$}
\put(7.3,.5){$\ol a_{n-2}$}
\put(8,1){$\bullet$}
\end{picture}
\end{center}

subject to the relations
$$a_0\ol a_0=\ol a_{n-2}a_{n-2}=\ol a_2a_2=0,\;\;
\ol a_ia_i=a_{i+1}\ol a_{i+1}\forall\;i\in\{5,\dots,n-2\}
$$
$$
\ol a_1a_1+a_2\ol a_2+a_3\ol a_3+f(\ol a_1a_1,a_2\ol a_2)=(\ol a_1a_1+a_2\ol a_2)^{n-3}=0
$$
for some
$$f\in rad^2(K<X,Y>/(X^3,Y^2,(X+Y)^{n-3})$$
so that
$$(X+Y+f(X,Y))^{n-3}=0.$$

\medskip

For all deformed preprojective algebras we number the
vertices by the condition that the vertex
$a_i$ starts at $v_i$ and ends at a vertex of higher label. This convention
numbers the vertices in a unique way.

\begin{Rem} \label{BESmistakes}
Observe that in \cite[page 238]{periodic} for type $E$ only
''admissible deformations'' may be applied,
which is the condition that  $(X+Y+f(X,Y))^{n-3}=0.$
However one relation is missing in \cite{periodic}  for type $D$ and type $E$
whereas the relation is correctly
displayed in \cite{BES}. I am grateful to Karin Erdmann for a clarifying email
on this subject.
\end{Rem}

The result is the following.

\begin{Theorem} (Bia\l kowski, Erdmann and Skowro\'nski \cite[Theorem 1.2]{BES})
Let $A$ be a finite dimensional selfinjective $K$-algebra.
Then $\Omega^3(S)\simeq S$ for every simple $A$-module if and only if
$A$ is preprojective of type $A_n$ for $(n\geq 1)$ or deformed preprojective of type
$D_n^f$ for $(n\geq 4)$, $E_6^f$, $E_7^f$, $E_8^f$ or $L_n^p$ for $(n\geq 1)$.
\end{Theorem}

It is a non trivial task to determine when deformations $f$ actually lead to non isomorphic
algebras. In a lecture at the ICRA XIV conference in Tokyo
in August 2010 Bia\l kowski annonce
that in characteristic $2$ the
deformed preprojective algebras of type $L_n^{X^{2j}}$ for $j\in\{0,1,\dots,n-1\}$
form a complete set of Morita equivalence classes
of these algebras, a fact that Skowro\'nski pointed out in an email to the
author from March 2007. Skowro\'nski announces in another email to the author in October 2008
that all algebras $L_n^p$ are symmetric and that moreover,
in characteristic different from $2$
the algebra $L_n^p$ is Morita equivalent to $L_n^{X^{n-1}}$, the non deformed preprojective algebra of type $L_n$.
The content of
Bia\l kowski's ICRA lecture are available in the conference abstract volume.

For type $D_n^f$ Bia\l kowski, Erdmann and Skowro\'nski
\cite[Proposition 6.2]{BES} show that the algebras
$D_n^{(XY)^j}$ are not Morita equivalent for different values of $j$.

No statement is known for type $E$ preprojective algebras.

\medskip

In joint work with Holm we computed the K\"ulshammer ideals for the algebras $L_n^{X^{2j}}$.
One main difficulty was to determine the commutator subspace. It is not very difficult to
get a generating set for the quotient of the algebra modulo the commutator space,
it is much more complicated to prove that the commutators one found really generate the
commutator space. In order to do so we apply the method described in
Section~\ref{Nakayamatwistedcentresection} and Section~\ref{Nakayamainpractice}.

A first step is

\begin{Lemma}\label{basisofcommutatorquotient} \cite[Lemma 3.9]{preproj}
Let $K$ be any field. Then
$L_{n+1}^p/[L_{n+1}^p,L_{n+1}^p]$ has a $K$-linear generating set
$$\{e_0,e_1,\dots,e_n,\epsilon, \epsilon^3,\epsilon^5,\epsilon^7,\dots, \epsilon^{2n+1}\}.$$
\end{Lemma}

We first need to fix a basis.

\begin{Prop}\label{abasis} \cite[Proposition 3.1]{preproj}
Let $K$ be any field.
A $K$-basis of $L_n^p$ is given by the following paths
between the vertices $i$ and $j$, where $i,j\in\{0,1,\ldots,n-1\}$.
$$\begin{array}{ll}
(1)\mbox{~~}a_ia_{i+1}\dots a_{j-1}
& \mbox{for~~} i<j \\
(2)\mbox{~~}a_ia_{i+1}\dots a_{j-1}a_j\dots
a_\ell\ol a_\ell \ol a_{\ell-1}\dots \ol a_j
& \mbox{for~~} i<j\mbox{~~and some~~}j\le \ell\le n-2 \\
(3)\mbox{~~}\ol a_{i-1}\ol a_{i-2}\dots \ol a_{j}
& \mbox{for~~}i\ge j \\
(4)\mbox{~~}a_ia_{i+1}\dots  a_\ell\ol a_\ell \ol a_{\ell-1}\dots
\ol a_i\dots\ol a_j
& \mbox{for~~}i \ge j\mbox{~~and some~~}i\le \ell\le n-2 \\
(5)\mbox{~~}
\ol a_{i-1}\ol a_{i-2}\dots \ol a_{0}\epsilon a_0 a_1\dots a_{j-1}
& \mbox{for any~~}i,j \\
(6)\mbox{~~}
\ol a_{i-1}\ol a_{i-2}\dots \ol a_{0}\epsilon a_0 a_1\dots
a_{\ell-1}a_{\ell}\ol a_{\ell}\ol a_{\ell-1}\dots \ol a_{j}
& \mbox{for~~} i<j\mbox{~~and some~~}j\le \ell \le n-2 \\
(7)\mbox{~~}
a_ia_{i+1}\dots a_\ell\ol a_\ell \ol a_{\ell-1}\dots \ol a_1\ol a_0\epsilon a_0a_1\dots a_{j-1}
& \mbox{for~~}i\ge j\mbox{~~and some~~}i\le \ell\le n-2
\end{array}
$$
\end{Prop}

Now we may define a Frobenius form with respect to this basis using Proposition~\ref{prop:form}.
It turns out that this form is in fact symmetric, non degenerate and associative.

\begin{Theorem} \cite[Theorem 3.5]{preproj}
\label{Lnissymmetric}
Let $K$ be any field.
Then the algebra $L_n^p$ is symmetric.
\end{Theorem}

We get many elements in the centre of $L_n^p$ by

\begin{Lemma}\label{centreLn} \cite[Lemma 3.13]{preproj}
Let $K$ be any field. Then
$$Z(L_n^p)\ni \epsilon^2+\epsilon^{3}p(\epsilon)+\sum_{j=0}^{n-2}(-1)^{j+1}\ol a_ja_j.$$ Hence
$$\left\{\left(\epsilon^2+\epsilon^{3}p(\epsilon)+
\sum_{j=0}^{n-2}(-1)^{j+1}\ol a_ja_j\right)^\ell\;|\;
\ell\in\{0,1,\dots,n-1\}\right\}\subseteq Z(L_n^p)$$ is a $K$-free subset.
Moreover $soc(Z(L_n^p))\subseteq Z(L_n^p)$.
\end{Lemma}

Lemma~\ref{centreLn}
provides a large space in $Z(L_n^p)$ of dimension
$2n$. Lemma~\ref{basisofcommutatorquotient}
shows that $L_n^p/[L_n^p,L_n^p]$ has dimension at most $2n$, whereas
Proposition~\ref{propertiesofNakayamatwistedcentre} show that the two
vector spaces are isomorphic. Hence the elements displayed in
Lemma~\ref{basisofcommutatorquotient} form a basis of $L_n^p/[L_n^p,L_n^p]$.

Using this result it is then possible to show

\begin{Theorem} \label{mainresult} (Holm and Zimmermann \cite[Theorem 4.1]{preproj})
Let $K$ be a perfect field of characteristic $2$.
Then for $0\leq j< n$ we get
\begin{eqnarray*}
\dim\left(T_i(L_n^{X^{2j}})\right)-\dim\left([L_n^{X^{2j}},L_n^{X^{2j}}]\right)
&=&
n-\max\left(\left\lceil\frac{2n-(2^{i+1}-2)j-(2^{i+1}-1)}{2^{i+1}}\right\rceil,\;\;0\right)
\end{eqnarray*}
\end{Theorem}

The attentive reader remarks that for the theorem one needs to suppose that the field
$K$ is perfect, whereas this was not supposed in the auxiliary steps. The reason for this
assumption comes from the technicalities in the proof of Theorem~\ref{mainresult}.
We need to find elements $x$ so that $x^{2^n}\in [L_n^p,L_n^p]$.
By what preceded and rather easy arguments it is necessary to consider this question only
for $x$ being a $K$-linear combination of odd powers of $\epsilon$. Now,
the $2^n$ powers of $x$ will give an element which is expressed in $2^n$-th powers
of the original coefficients, and from there it is not hard to imagine that one needs to
take $2^n$-th roots of the solutions, in order to get the original coefficients from some expression
one obtains from some solution one got by a linear algebra argument.

Moreover, we remark that the Morita equivalence classification of the algebras $L_n^p$ is finer
than what can be done by K\"ulshammer ideals. However, K\"ulshammer ideals
distinguish algebras up to derived equivalences and
even up to stable equivalences of Morita type (cf Theorem~\ref{LZZtheorem}).

\subsection{Algebras of polynomial growth and domestic weakly symmetric algebras}

\label{polygrowthsection}

Let $K$ be an algebraically closed field and let $A$ be a finite dimensional
$K$-algebra.
\begin{itemize}
\item
The algebra $A$ is called of {\em finite representation type} if
$A$ admits only a finite number of indecomposable $A$-modules up to isomorphism.
\item
The algebra $A$ is called of {\em tame representation type} if $A$
is not of finite representation type and if for every positive
integer $d$ there are a finite number of
$A\otimes_KK[X]$-modules $M_1(d),M_2(d),\dots,M_{n(d)}(d)$,
which are free as $K[X]$-modules and so that
for each $d$
all but a finite number of $d$-dimensional indecomposable $A$-modules are isomorphic
to a module of the form $M_i(d)\otimes_KK[X]/(X-\lambda)$ for some $\lambda\in K$ and
some $i\in\{1,\dots,n(d)\}$.
\begin{itemize}
\item
The tame algebra $A$ is called of {\em domestic representation type} if,
taking $n(d)$ as small as possible, there is an integer $m$ so that $n(d)\leq m$
for all $d$.
\item
The tame algebra $A$ is called of {\em polynomial growth} if, taking
$n(d)$ as small as possible, there is an integer $m$
so that $$\lim_{n\ra\infty}\frac{n(d)}{d^m}=0.$$
\end{itemize}
\item
The algebra $A$ is called of {\em wild representation type} if for every algebra $B$
there is a functor $B-mod\lra A-mod$ which is exact, preserves isomorphism classes
and carries indecomposable objects to indecomposable objects.
\end{itemize}

Of course, tame domestic algebras are polynomial growth tame algebras.

A fundamental result of Drozd says that $A$ is either tame or wild or of finite
representation type. There is intensive research aiming a possible classification
of algebras of tame representation type, though the goal seems still to be very far.
Nevertheless Bocian, Holm and Skowro\'nski classified tame domestic weakly symmetric
algebras \cite{BocianHolmSkowronski2004,BocianHolmSkowronski2007,weaklyprojdomestic}
and tame weakly symmetric algebras of polynomial growth \cite{BialkowskiHolmSkowronski2003a,BialkowskiHolmSkowronski2003b,weklyprojpolygrowth}
up to derived equivalences in a series of papers.

We present some of the details.

\begin{Def}
A selfinjective algebra of tame representation type is called standard if
its basic algebra
admits simply connected Galois coverings.
Else a selfinjective algebra is called non standard.
\end{Def}

\begin{Theorem}\cite[Theorem 1]{BocianHolmSkowronski2004}
\label{DomesticStandardSingular}
For an algebra $A$ the following
statements are equivalent:
\begin{itemize}
\item[(1)] $A$ is a representation-infinite domestic selfinjective
algebra having simply connected Galois coverings and the Cartan
matrix $C_A$ is singular.

\item[(2)] $A$ is derived equivalent to the trivial extension
$T(C)$ of a canonical algebra $C$ of Euclidean type.

\item[(3)] $A$ is stably equivalent to the trivial extension
$T(C)$ of a canonical algebra $C$ of Euclidean type.
\end{itemize}
Moreover, the trivial extensions $T(C)$ and $T(C')$ of two
canonical algebras $C$ and $C'$ of Euclidean type are derived
equivalent (respectively, stably equivalent) if and only if the
algebras $C$ and $C'$ are isomorphic.
\end{Theorem}

In order to be able to formulate  Bocian, Holm and Skowro\'nski's result
for weakly symmetric algebras of domestic
representation type with nonsingular Cartan matrices,  we need to
define the following algebras.

$$\unitlength0.6cm
\begin{picture}(11,6)
 \put(9.4,3){\circle{2.0}}
\put(6.8,3){\circle{2.0}} \put(7.85,2.9){\vector(0,1){0.3}}
\put(8,2.9){$\bullet$} \put(8.35,2.9){\vector(0,1){0.3}}
\put(4.8,2.9){$\beta$} \put(11,2.9){$\alpha$}

 \put(0,1.9){$A(\lambda)$}
\put(0,1){$\lambda\in K\backslash \{0\}$}

\put(4.8,0.5){$\alpha^2=0, \beta^2=0,
\alpha\beta=\lambda\beta\alpha$}
\end{picture}$$

$$\unitlength0.5cm
\begin{picture}(30,15)
 \put(18.1, 11.1){\vector(1,1){1.8}}

\put(16.5, 11.5){$\beta_1$}

 \put(18, 11){\vector(-1,1){1.8}}

 \put(14.5, 13){$\beta_2$}

\put(16, 13){\vector(-2,1){1.8}}

\put(13, 13.2){$\beta_3$}

\put(14, 14){\vector(-1,0){1.8}}

\put(11, 13.0){$\beta_4$}

\put(12, 14){\vector(-2,-1){1.8}}

\put(10.7, 8.8){$\beta_{q-3}$}

\put(10, 9){\vector(2,-1){1.8}}

\put(12.2, 8.5){$\beta_{q-2}$}

\put(12, 8){\vector(1,0){1.8}}

\put(14, 8.8){$\beta_{q-1}$}

\put(14, 8){\vector(2,1){1.8}}

\put(16.4, 10.3){$\beta_{q}$}

\put(16, 9){\vector(1,1){1.8}}

\put(19.2, 11.5){$\alpha_{1}$}

\put(20, 9){\vector(-1,1){1.8}}

\put(21, 13){$\alpha_{2}$}

\put(22, 8){\vector(-2,1){1.8}}

\put(22.5, 13.5){$\alpha_{3}$}

\put(24, 8){\vector(-1,0){1.8}}

\put(24.5, 13){$\alpha_{4}$}

\put(26, 9){\vector(-2,-1){1.8}}

\put(18.8, 10.3){$\alpha_{p}$}

\put(20.7, 8.8){$\alpha_{p-1}$}

\put(22.4, 8.5){$\alpha_{p-2}$}

\put(24.1, 8.8){$\alpha_{p-3}$}

\put(20, 13){\vector(2,1){1.8}}

\put(22, 14){\vector(1,0){1.8}}

\put(24, 14){\vector(2,-1){1.8}}

\multiput(9, 11)(0.1,0.2){10}{\circle*{0.01}}

\multiput(9, 11)(0.1,-0.2){10}{\circle*{0.01}}

\multiput(27, 11)(-0.1,0.2){10}{\circle*{0.01}}

\multiput(27, 11)(-0.1,-0.2){10}{\circle*{0.01}}

 \put(1, 12){$A(p, q)$}

 \put(1, 10.5){$1\leq p\leq q$}
 \put(1, 9){$p+q\geq 3$}

 \put(9, 6){$\alpha_1\alpha_2\cdots \alpha_p\beta_1\beta_2\cdots \beta_q=\beta_1\beta_2\cdots \beta_q\alpha_1\alpha_2\cdots \alpha_p$}

\put(9, 4.5){$ \alpha_p\alpha_1=0, \beta_q\beta_1=0$}

\put(9, 3){$\alpha_i\alpha_{i+1}\cdots \alpha_p\beta_1\cdots
\beta_q  \alpha_1 \cdots \alpha_{i-1}\alpha_i=0, 2\leq i\leq p$}

\put(9, 1.5){$\beta_j\beta_{j+1}\cdots \beta_q\alpha_1\cdots
\alpha_p \alpha_1 \cdots \beta_{i-1}\beta_i=0, 2\leq j\leq q$}

\end{picture}$$

$$\unitlength0.5cm
\begin{picture}(30,11)

\multiput(9, 7)(0.1,0.2){10}{\circle*{0.01}}

\multiput(9, 7)(0.1,-0.2){10}{\circle*{0.01}}

\put(16.5, 7.5){$\beta_1$}

 \put(18, 7){\vector(-1,1){1.8}}

 \put(14.5, 9){$\beta_2$}

\put(16, 9){\vector(-2,1){1.8}}

\put(13, 9.2){$\beta_3$}

\put(14, 10){\vector(-1,0){1.8}}

\put(11, 9.0){$\beta_4$}

\put(12, 10){\vector(-2,-1){1.8}}

\put(10.7, 4.8){$\beta_{n-3}$}

\put(10, 5){\vector(2,-1){1.8}}

\put(12.2, 4.5){$\beta_{n-2}$}

\put(12, 4){\vector(1,0){1.8}}

\put(14, 4.8){$\beta_{n-1}$}

\put(14, 4){\vector(2,1){1.8}}

\put(16.4, 6.3){$\beta_{n}$}

\put(16, 5){\vector(1,1){1.8}}

\put(19.2, 7){\circle{2}}

\put(19, 8){\vector(-1,0){0.01}}

\put(20.5, 7){$\alpha$}

 \put(1, 8){$\Lambda(n)$}

 \put(1, 6.5){$n\geq 2$}

 \put(9, 2){$\alpha^2=(\beta_1\beta_2\cdots\beta_n)^2, \alpha\beta_1=0, \beta_n\alpha=0$}

\put(9, 0.5){$\beta_j\beta_{j+1}\cdots \beta_n\beta_1\cdots
\beta_n \alpha_1 \cdots \beta_{i-1}\beta_i=0, 2\leq j\leq n$}

\end{picture}$$

$$\unitlength0.5cm
\begin{picture}(30,15)

\multiput(9, 11)(0.1,0.2){10}{\circle*{0.01}}

\multiput(9, 11)(0.1,-0.2){10}{\circle*{0.01}}

\put(16.5, 11.5){$\beta_1$}

 \put(18, 11){\vector(-1,1){1.8}}

 \put(14.5, 13){$\beta_2$}

\put(16, 13){\vector(-2,1){1.8}}

\put(13, 13.2){$\beta_3$}

\put(14, 14){\vector(-1,0){1.8}}

\put(11, 13.0){$\beta_4$}

\put(12, 14){\vector(-2,-1){1.8}}

\put(10.7, 8.8){$\beta_{n-3}$}

\put(10, 9){\vector(2,-1){1.8}}

\put(12.2, 8.5){$\beta_{n-2}$}

\put(12, 8){\vector(1,0){1.8}}

\put(14, 8.8){$\beta_{n-1}$}

\put(14, 8){\vector(2,1){1.8}}

\put(16.4, 10.3){$\beta_{n}$}

\put(16, 9){\vector(1,1){1.8}}

\put(18.2, 11.3){\vector(1,2){1.2}}

\put(19.8, 13.4){\vector(-1,-2){1.2}}

\put(18.4, 10.9){\vector(1,-2){1.2}}

\put(19.2, 8.4){\vector(-1,2){1.2}}

\put(19.4, 9.7){$\gamma_1$}

\put(19.4, 11.7){$\alpha_2$}

\put(18, 9){$\gamma_2$}

\put(18, 124.5){$\alpha_1$}

 \put(1, 12){$\Gamma(n)$}

 \put(1, 10.5){$n\geq 1$}

 \put(9, 6){$\alpha_1\alpha_2=(\beta_1\beta_2\cdots\beta_n)^2=\gamma_1\gamma_2,$}

\put(9, 4.5){$\alpha_2\beta_1=0, \gamma_2\beta_1=0,
\beta_n\alpha_1=0$}

\put(9, 3){$\beta_n\gamma_1=0, \alpha_2\gamma_1=0,
\gamma_2\alpha_1=0$}

\put(9, 1.5){$\beta_j\beta_{j+1}\cdots \beta_n\beta_1\cdots
\beta_n \alpha_1 \cdots \beta_{i-1}\beta_i=0, 2\leq j\leq n$}

\end{picture}$$

\begin{Theorem} (Bocian, Holm and Skowro\'nski \cite[Theorem
2]{BocianHolmSkowronski2004})\label{DomesticStandardNonSingular}
For a domestic standard selfinjective algebra $A$   the
following statements are equivalent:
\begin{itemize}
\item[(1)]$A$ is weakly symmetric and the Cartan matrix $C_A$ is
nonsingular.

\item[(2)] $A$ is derived equivalent to an algebra of the form
$A(\lambda), A(p, q), \Lambda(n),  \Gamma(n)$.

\item[(3)] $A$ is stably equivalent to an algebra of the form
$A(\lambda), A(p, q), \Lambda(n), \Gamma(n)$.

\end{itemize}

Moreover, two algebras of the forms $A(\lambda), A(p, q),
\Lambda(n), \Gamma(n)$ are derived equivalent (respectively,
stably equivalent) if and only if they are isomorphic.
\end{Theorem}

For algebras of non standard type we need to introduce the
following algebra $\Omega(n)$.
The quiver with relations of $\Omega(n)$ is as
follows.

$$\unitlength0.5cm
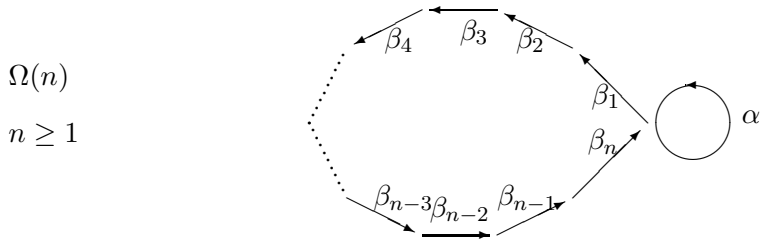
\begin{picture}(30,11)

\multiput(9, 7)(0.1,0.2){10}{\circle*{0.01}}

\multiput(9, 7)(0.1,-0.2){10}{\circle*{0.01}}

\put(16.5, 7.5){$\beta_1$}

 \put(18, 7){\vector(-1,1){1.8}}

 \put(14.5, 9){$\beta_2$}

\put(16, 9){\vector(-2,1){1.8}}

\put(13, 9.2){$\beta_3$}

\put(14, 10){\vector(-1,0){1.8}}

\put(11, 9.0){$\beta_4$}

\put(12, 10){\vector(-2,-1){1.8}}

\put(10.7, 4.8){$\beta_{n-3}$}

\put(10, 5){\vector(2,-1){1.8}}

\put(12.2, 4.5){$\beta_{n-2}$}

\put(12, 4){\vector(1,0){1.8}}

\put(14, 4.8){$\beta_{n-1}$}

\put(14, 4){\vector(2,1){1.8}}

\put(16.4, 6.3){$\beta_{n}$}

\put(16, 5){\vector(1,1){1.8}}

\put(19.2, 7){\circle{2}}

\put(19, 8){\vector(-1,0){0.01}}

\put(20.5, 7){$\alpha$}

 \put(1, 8){$\Omega(n)$}

 \put(1, 6.5){$n\geq 1$}

 \put(3, 2){$\alpha^2=\alpha\beta_1\beta_2\cdots\beta_n,
 \alpha\beta_1\beta_2\cdots\beta_n+\beta_1\beta_2\cdots\beta_n\alpha=0, $}

\put(3, 0.5){$\beta_n\beta_{1}=0, \beta_j\beta_{j+1}\cdots
\beta_n\beta_1\cdots \beta_n \alpha_1 \cdots \beta_{i-1}\beta_i=0,
2\leq j\leq n$}

\end{picture}$$

\begin{Theorem} (Bocian, Holm and Skowro\'nski
\cite[Theorem 1]{BocianHolmSkowronski2007})\label{DomesticNonStandard}
Any nonstandard representation infinite
selfinjective algebra of domestic type is derived equivalent
(resp. stably equivalent) to an algebra $\Omega(n)$ with $n\geq
1$. Moreover, two algebras $\Omega(n)$ and $\Omega(m)$ are derived
equivalent (respectively, stably equivalent) if and only if $m =n$.
\end{Theorem}

Bocian, Holm and Skowro\'nski show that standard and non standard
domestic algebras cannot be derived equivalent. We are able
to improve partially the result.

\begin{Lemma} (Zhou and Zimmermann \cite[Lemma 2.3]{polygrowth})
\label{WeaklySymmetricDomesticStandardVsNonStandard}
 A standard weakly symmetric algebra of domestic type
cannot be stably equivalent to a nonstandard one.
\end{Lemma}

The method of proof is to compare stable Auslander Reiten quivers and
to use the fact that the class of special biserial algebras
is closed under stable equivalences, a result due to Pogorza\l y
\cite[Theorem 7.3]{Pogorzaly}.

Observe that a Morita equivalence classification of standard self-injective
domestic algebras is not given in the results of Bocian, Holm and Skowro\'nski.
This is the reason why we suppose that the domestic algebras
are weakly symmetric and not only selfinjective.

\bigskip

We show

\begin{Theorem} (Zhou and Zimmermann \cite[Theorem 2.5]{polygrowth})
\label{ClassificationSymmetricDomestic}
\begin{enumerate}
\item
Two weakly symmetric  algebras of domestic representation type   are
derived equivalent if and only
if they are stably equivalent.

\item The class of weakly symmetric
algebras  of domestic representation type  is closed under stable
equivalences.
\end{enumerate}
\end{Theorem}

Bocian, Holm and Skowro\'nski give a classification of symmetric algebras
of polynomial growth up to derived equivalences
\cite{BialkowskiHolmSkowronski2003a},
\cite{BialkowskiHolmSkowronski2003b},
\cite{weklyprojpolygrowth}. They get a finite list of algebras which are
defined by quivers and relations, and where the relations involve some parameters.
We are not completely able to distinguish the algebras with the same
quiver and relations and different parameters. We call this problem the scalar problem.

As for symmetric algebras of polynomial growth we get the following result.

\begin{Theorem} (Zhou and Zimmermann \cite[Theorem 3.5]{polygrowth}) \label{NonDomesticClassification}
The classification  of  indecomposable
non-domestic weakly symmetric algebras of polynomial growth
up to stable equivalences of Morita type coincide with the derived
equivalence classification, modulo the scalar problems.
\end{Theorem}

The method of proof uses, among other general arguments, a
computation of K\"ulshammer ideals.

Concerning the Auslander Reiten conjecture we get

\begin{Cor}
(Zhou and Zimmermann \cite[Corollary 2.7, Theorem 3.6, Theorem 3.7]{polygrowth})
\begin{itemize}
\item Let $A$ be an indecomposable algebra stably equivalent
  to an indecomposable   symmetric algebra $B$ of
domestic type. Then $A$ and $B$ have the same number of isomorphism classes of simple modules.
\item Let $A$ be an indecomposable algebra stably equivalent
  to an indecomposable weakly symmetric standard algebra $B$ of
domestic type. Then $A$ and $B$ have the same number of isomorphism classes of simple modules.
\item
Let $A$ and $B$ be indecomposable algebras which are stably equivalent of Morita type.
If $A$ is tame symmetric with only $\Omega$-periodic modules, then $A$ and $B$
have the same number of isomorphism classes of simple modules.
\item
Let $A$ be an indecomposable
algebra and suppose $A$ is stably equivalent of Morita type to an indecomposable
non-domestic  symmetric algebra  $B$   of polynomial growth.
Then $A$ and $B$ have the same number of isomorphism classes of simple modules.
\item
Let $A$ and  $B$  be two indecomposable algebras which
are standard non-domestic  weakly symmetric algebra of polynomial
growth or non-standard non-domestic  selfinjective algebra of
polynomial growth. If they are stably equivalent of Morita type,
then $A$ and $B$ have the same number of isomorphism classes of simple modules.
\end{itemize}
\end{Cor}

\section{Hochschild (co-)homology}

\label{Hochschildsection}

It is well-known (cf e.g. Loday \cite[Paragraphs 1.1.6 and 1.5.2]{Loday}),
and actually we already used this fact implicitly
in e.g. in Section~\ref{Nakayamatwistedcentresection} and
Section~\ref{stableinvariance}, that
$A/[A,A]=HH_0(A)$ is the degree $0$ Hochschild homology
and $Z(A)=HH^0(A)$ is the degree $0$ Hochschild cohomology.
A natural question becomes now to generalise K\"ulshammer's ideals to higher Hochschild
(co-)homology. This was done by the author in \cite{Gerstenhaber} and \cite{TAHochschild}.

The symmetrising form on $A$ induces a non degenerate pairing $HH^0(A)\times HH_0(A)\lra K$.
Hence, in order to generalise to higher Hochschild (co-)homology  we first
need to produce a non degenerate bilinear form
$$HH^n(A)\times HH_n(A)\lra K$$
for symmetric algebras induced by the symmetrising form.
Let $\B A$ be the bar resolution (i.e. a specific projective resolution; and
actually projective resolution will do at this place) of $A$ as
$A\otimes_KA^{op}$-modules.

\begin{eqnarray*}
Hom_K(HH_n(A),K)&=&Hom_K(H_n(\B A\otimes_{A\otimes A^{op}}A),K)\\
&=&H^n(Hom_K(\B A\otimes_{A\otimes A^{op}}A,K))\\
&=&H^n(Hom_{A\otimes_KA^{op}}(\B A,Hom_K(A,K))\\
&\simeq&H^n(Hom_{A\otimes_KA^{op}}(\B A,A))\\
&=&HH^n(A)
\end{eqnarray*}
where the second last isomorphism is induced by the symmetrising form $A\simeq Hom_K(A,K)$
as $A\otimes_KA^{op}$ modules.

This isomorphism yields a non degenerate bilinear form
$$\langle\;,\;\rangle_n:HH^n(A)\times HH_n(A)\lra K$$
which extends the symmetrising form on $A$.

In order to define K\"ulshammer ideals we used the $p$-power map
on $HH_0(A)=A/[A,A]$. However, how to get a $p$-power map on $HH_*(A)$ is not
completely clear.
Nevertheless, the multiplicative structure on $HH^*(A)$
defined by the cup product can be used instead. By adjointness
with respect to the bilinear form we get then an analogue
of the K\"ulshammer ideal structure on Hochschild homology
instead of cohomology.
This dual construction on $A/[A,A]=HH_0(A)$ was studied by
K\"ulshammer as well in \cite[Part IV]{Ku1}.

For the $p^n$-power mapping by cup product
$HH^{m}(A,A)\lra HH^{p^nm}(A,A)$ one gets a right adjoint
$\kappa_n^{(m),A}:HH_{p^nm}(A,A)\lra HH_{m}(A,A)$ with respect
to $\langle\;,\;\rangle_m$ and $\langle\;,\;\rangle_{p^nm}$
Observe that, if $p$ is odd, the cup product square is $0$ in
odd degree cohomology.
Hence for $p$ odd the $p$-power map as well as the adjoint is $0$
in odd degree cohomology.

\begin{Prop} \cite[Lemma 2.6, beginning remarks of Section 3]{Gerstenhaber}
Let $A$ be a finite dimensional symmetric $K$ algebra over
a perfect field $K$ of characteristic $p>0$. Then,
for all ${n\in\N}$ and for all ${x\in HH_{p^nm}(A,A)}$ there is a unique
$\kappa_n^{(m)}(x)\in HH_m(A,A)$ so that for all $f\in HH^m(A,A)$
one has
$$\left\langle f^{p^n},x\right\rangle_{p^nm}=
\left(\left\langle f,\kappa_n^{(m)}(x)\right\rangle_m\;\right)^{p^n}.$$
\end{Prop}

Using the mapping $\kappa_n^{(m)}$ one gets the compatibility
with derived equivalences.

\begin{Theorem}\label{derivedkappacup}\cite[Theorem 1]{Gerstenhaber}
Let $A$ be a finite dimensional symmetric $K$-algebra over a perfect field
$K$ of characteristic $p>0$. Let $B$ be a second algebra so that
$D^b(A)\simeq D^b(B)$ as triangulated categories.
Let $p$ be a prime and let $m\in\N$.
Then, there is a standard equivalence $F:D^b(A)\simeq D^b(B)$,
and any such standard equivalence induces an isomorphism
$HH_m(F):HH_m(A,A)\lra HH_m(B,B)$ of all Hochschild homology
groups satisfying
$$HH_{m}(F)\circ\kappa_n^{(m),A}\circ HH_{p^nm}(F)^{-1}=
\kappa_n^{(m),B}\;.$$
\end{Theorem}

The proof is very much like the one in degree $0$ for Hochschild cohomology.
Nevertheless a clear definition of an isomorphism on Hochschild homology
induced by a standard derived equivalence was not published explicitly before.
The construction was somewhat implicit in Rickard's work, but an explicit
clarification seems to appear in \cite[Section 1.2]{Gerstenhaber} for the first time.

One should notice however that a derived equivalence may be non standard,
and then it is not clear how to define an induced mapping on the category
of bimodules, and in the sequel on the Hochschild homology. One needs the
standard equivalence in order to control the way it acts on Hochschild homology.

What happens if $A$ is not symmetric? Already for the non degenerate pairing
between Hochschild homology and cohomology it is
not clear how to to define it properly. We may again use trivial extension algebras.
Then there are ring homomorphisms $A\lra \T A$ and $\T A\lra A$
by projection, injection from and to the second component. Hochschild homology
is functorial, contrary to Hochschild cohomology. Hence we get
mappings
$$HH_n(\iota_A):HH_n(A)\lra HH_n(\T A)$$
and
$$HH_n(\pi_A):HH_n(\T A)\lra HH_n(A).$$

Now, defining
$$\hat\kappa_n^{(m)}:=HH_m(\pi_A)\circ \kappa_n^{(m)}\circ HH_{p^nm}(\iota_A):
HH_{p^nm}(A)\lra HH_m(A)$$
one obtains an invariant under derived equivalences.

\begin{Theorem}\cite[Theorem 2]{TAHochschild}\label{TAHochschildthm}
Let $K$ be a perfect field of characteristic $p>0$, let
$A$ and $B$ be finite dimensional  $K$-algebras and suppose
that $D^b(A)\simeq D^b(B)$ as triangulated categories.
Let $F$ be an explicit standard equivalence between $D^b(A)$ and
$D^b(B)$. Then, $F$ induces a sequence of isomorphisms
$HH_m(F):HH_m(A)\lra HH_m(B)$ so that
$$HH_m(F)\circ \hat \kappa_n^{(m);A}=\hat \kappa_n^{(m);B}\circ HH_{p^nm}(F).$$
\end{Theorem}

Obviously Theorem~\ref{TAHochschildthm} generalises
Theorem~\ref{derivedkappacup} to non symmetric algebras.
Since Hoch\-schild homology is often better understood than Hochschild
cohomology we expect that this generalisation will bear use in the future.

Nevertheless there is a $p$-power map available in some cases coming from
the Gerstenhaber structure on Hochschild homology. The construction is due to
Stasheff and Quillen.

Let
$$Coder(\B(A),\B(A)):=
\{D\in End_{A\otimes A^{op}}(\B(A))| \Delta\circ D=
(id_{\B(A)}\otimes D+D\otimes id_{\B(A)})\circ \Delta\}$$
be the coderivations. Since $\B(A)$ is graded, $Coder(\B(A),\B(A))$ is graded
as well. Denote by $Coder^n(\B(A),\B(A))$  the degree $n$ coderivations.
The vector space $Coder(\B(A),\B(A))$ is a
graded Lie algebra with Lie bracket being the commutator.
Moreover (cf  Stasheff \cite[Proposition]{Stasheff}),
$$Coder(\B(A),\B(A))\simeq Hom_{A\otimes A^{op}}(\B(A),A)[1]\;.$$

The key observation is

\begin{Lemma} \cite[Lemma 4.1]{Gerstenhaber}\label{Gerstenhaberpower}
(Keller, personal communication)
\begin{itemize}
\item
Suppose $K$ is a field. Then,
$$D\in Coder^{2n+1}(\B(A),\B(A))\Rightarrow
D^{2}\in Coder^{2\cdot(2n+1)}(\B(A),\B(A)).$$
\item
Suppose $K$ is a field of characteristic $p>0$.
Then, $$D\in Coder^{2n}(\B(A),\B(A))\Rightarrow
D^{p}\in Coder^{2pn}(\B(A),\B(A)).$$
\end{itemize}
\end{Lemma}

This $p$-power structure carries over to Hochschild cohomology.

\begin{Lemma}\label{gehtnicht}\cite[Lemma 4.2]{Gerstenhaber}
Let $K$ be a field of characteristic $p>0$ and let $D\in Coder^n(\B(A),\B(A))$.
\begin{enumerate}
\item If $p=2$ and $n\in\N$,
then the mapping $D\mapsto D^2$ induces a mapping
$$HH^{n+1}(A,A)\lra HH^{2n+1}(A,A)$$
\item
If $p>2$ and $n=2m\in 2\N$, then
the mapping $D\mapsto D^p$ induces a mapping
$$HH^{2m+1}(A,A)\lra HH^{2pm+1}(A,A)$$
\end{enumerate}
\end{Lemma}

Hence, for $p=2$ the Hochschild cohomology becomes a $2$-restricted Lie algebra
with the Gerstenhaber bracket and the $2$-power map. For $p>2$ the
odd degree Hochschild cohomology becomes a $p$-restricted Lie algebra
with the Gerstenhaber bracket and the $p$-power mapping.

Using these constructions we get

\begin{Theorem}\label{Gerstenhaberkompatibel}\cite[Proposition 4.4]{Gerstenhaber}
Let $A$ and $B$ be $K$-algebras over a field $K$.
Suppose $D^b(A)\simeq D^b(B)$
as triangulated categories.
\begin{itemize}
\item
If the characteristic of $K$ is $2$ then $HH^\ast(A,A)$ and
$HH^\ast(B,B)$ are isomorphic as restricted Lie algebras.
\item
If the characteristic of $K$ is $p>2$, then the Lie algebras
consisting of odd degree Hochschild cohomologies
$\bigoplus_{n\in\N}HH^{2n+1}(A,A)$ and $\bigoplus_{n\in\N}HH^{2n+1}(B,B)$
are isomorphic as restricted Lie algebras.
\end{itemize}
\end{Theorem}

However we fail to prove that the so-defined Gerstenhaber $p$-power map is
additive, neither semilinear. Hence it does not seem to be clear how a
K\"ulshammer structure could be built from there.

\begin{Acknowledgement}
I wish to express my thanks to Javad Asadollahi and all mathematicians at the
IPM in Teheran for their great hospitality during the conference
''Representation theory of Algebras'' in June 2008.

I also wish to thank Karin Erdmann and Andrzej Skowro\'nski for very useful
and clarifying discussions, and the referee for a very detailed, careful and
most helpful report.
\end{Acknowledgement}


\begin{thebibliography}{88}

\bibitem{AndreuEstefania}
Estefan\'\i a Andreu Juan, {\em The Hochschild cohomology ring
of preprojective algebras of type $L_n$}, preprint 2010; arXiv math.RT 1010.2790v1

\bibitem{AR}
Maurice Auslander and Idun Reiten, {\em
Stable equivalence of Artin algebras.}
Proceedings of the Conference on Orders, Group Rings and Related Topics
(Ohio State Univ., Columbus, Ohio, 1972), 8-71.
Springer Lecture Notes in Mathematics 353  (1973).

\bibitem{ARS}
Maurice Auslander, Idun Reiten and Sverre Sm\aa l\o{}, {\sc Representation Theory of
Artin Algebras}, Cambridge University Press 1995.

\bibitem{nonsymmetric}
Christine Bessenrodt, Thorsten Holm and Alexander Zimmermann,
{\em Generalized Reynolds ideals for non-symmetric algebras},
Journal of Algebra {\bf 312} (2007) 985-994.

\bibitem{BES}
Jerzy Bia\l kowski, Karin Erdmann and Andrzej Skowro\'nski, {\em Deformed
preprojective algebras of generalized Dynkin type},
Transactions of the American Mathematical Society, {\bf 359} (2007) 2625-2650.

\bibitem{BialkowskiHolmSkowronski2003a}
Jerzy  Bia\l kowski, Thorsten Holm and Andrzej Skowro\'nski,
\emph{Derived
equivalences for tame weakly symmetric algebras having only
periodic modules}. Journal of Algebra \textbf{269} (2003), no. 2, 652-668.

\bibitem{BialkowskiHolmSkowronski2003b}
Jerzy  Bia\l kowski, Thorsten Holm and Andrzej Skowronski,
\emph{On
nonstandard tame selfinjective algebras having only periodic
modules}. Colloquium Mathematicum \textbf{97} (2003), no. 1, 33-47.


\bibitem{Bleher}
Frauke Bleher, {\it Dihedral blocks with two simple modules},
Proceedings of the American Mathematical Society {\bf 138} (2010) 3467-3479.

\bibitem{BocianHolmSkowronski2004}
Rafa\l\ Bocian, Thorsten  Holm and  Andrzej Skowro\'nski,
{\em Derived equivalence classification of weakly symmetric algebras of Euclidean type}.
Journal of Pure and Applied Algebra  \textbf{191}  (2004),  no. 1-2, 43-74.


\bibitem{BocianHolmSkowronski2007}
Rafa\l\ Bocian, Thorsten  Holm and  Andrzej Skowro\'nski,
{\em  Derived equivalence classification of nonstandard selfinjective algebras of domestic type}.
Communications in Algebra  \textbf{35}  (2007),  no. 2, 515--526.

\bibitem{Bouc}
Serge Bouc, {\em Bimodules, trace g\'en\'eralis\'ee,
et transferts en homologie de Hochschild}, preprint available at
http://www.lamfa.u-picardie.fr/bouc/transfer.pdf

\bibitem{Brauer}
Richard Brauer, {\em Zur Darstellungstheorie der Gruppen endlicher Ordnung},
Mathematische Zeitschrift {\bf 63} (1956) 406-444.

\bibitem{BHHKM}
Thomas Breuer, L\'aszl\'o H\'ethelyi, Erzs\'ebet Horv\'ath, Burkhard K\"ulshammer,
John Murray, {\em Cartan invariants and central ideals of group algebras,}
Journal of Algebra {\bf 296} (2006) 177-195.

\bibitem{Broue1994}{{Michel Brou\'e,} {\em Equivalences of blocks of group
    algebras.} In: { Finite
dimensional algebras and related topics}. V.Dlab and L.L.Scott
(eds.), Kluwer, 1994, 1-26.}

\bibitem{CarlsonThevenaz}
Jon Carlson and Jacques Th\'evenaz, {\em The cassification of torsion endo-trivial modules},
Annals of Mathematics {\bf 162} (2005) 823-883

\bibitem{ErdmannLNM}
Karin Erdmann, {\sc Blocks of tame representation type and related
algebras}, Springer Lecture Notes in Mathematics 1428 (1990).

\bibitem{periodic}
Karin Erdmann and Andrzej Skowro\'nski, {\em Periodic algebras},
in ''Trends in representation theory of
algebras and related topics'', European Mathematical Society,
Series of Congress Reports, (2008) pp 201-251

\bibitem{Feit}
Walter Feit, {\sc The representation theory of finite groups}, North-Holland Amsterdam, 1982.

\bibitem{HHKM}
L\'aszl\'o H\'ethelyi, Ersz\'ebet Horv\'ath, Burkhard K\"ulshammer
and John Murray, {\em Central ideals and Cartan invariants of
symmetric algebras}, Journal of Algebra {\bf 293} (2005) 243-260.

\bibitem{Holmhabil}
Thorsten Holm, {\em Blocks of tame representation type and related
algebras: derived equivalences and Hochschild cohomology},
Habilitationsschrift (2001) Universit\"at Magdeburg.

\bibitem{weaklyprojdomestic}    Thorsten  Holm and  Andrzej Skowro\'nski,
{\em     Derived equivalence classification of symmetric algebras of domestic type}.
Journal of the Mathematical Society of Japan  \textbf{58}  (2006),  no. 4, 1133--1149


\bibitem{weklyprojpolygrowth}    Thorsten  Holm and  Andrzej Skowro\'nski,
{\em    Derived equivalence classification of symmetric algebras of  polynomial growth}.
to appear in Glasgow Journal of Mathematics.

\bibitem{hz-tame}
Thorsten Holm and Alexander Zimmermann,
{\em Generalized Reynolds ideals and  derived equivalences
for algebras of dihedral and semidihedral type},
Journal of Algebra {\bf 320} (2008)
3425-3437.

\bibitem{preproj} Thorsten Holm and Alexander Zimmermann,
{\em Deformed preprojective
algebras of type $L$: K\"ulshammer spaces and derived equivalences}, preprint 2010.


\bibitem{KR}
Michael Kauer, Klaus\,W. Roggenkamp, {\em Higher-dimensional orders, graph-orders,
and derived equivalences.} Journal of Pure and Applied Algebra \textbf{155} (2001),
181-202.

\bibitem{Keller}
Bernhard Keller, {\em Deriving DG categories}, Annales
Scientifiques Ecole Normale Sup\'erieure {\bf 27} (1994), 63 - 102.

\bibitem{KellerVossieck}  Bernhard Keller and Dieter Vossieck,
\emph{Sous les cat\'egories d\'eriv\'ees}.
Comptes Rendus de l'Acad\'emie des Sciences Paris, S\'erie I Math\'ematique
\textbf{305}  (1987),  no. 6, 225--228.

\bibitem{derbuch}
Steffen K\"onig and Alexander Zimmermann, {\sc Derived Equivalences for Group Rings},
with contributions by Bernhard Keller, Markus Linckelmann, Jeremy Rickard and Rapha\"el Rouquier; Springer Lecture Notes in Mathematics 1685 (1998).

\bibitem{KLZ}
Steffen K\"onig, Yuming Liu and Guodong Zhou, {\em Transfer maps
in Hochschild (co-)homology and applications to stable and derived
invariants and to the Auslander-Reiten conjecture, } to appear
in Transactions of the American Mathematical Society.

\bibitem{Ku1}
 Burkhard K\"ulshammer,
{\em Bemerkungen \"uber die Gruppenalgebra als symmetrische Algebra
I,\,II,\,III,\,IV},
Journal of Algebra \textbf{72} (1981), 1--7;
Journal of Algebra \textbf{75} (1982), 59--69;
Journal of Algebra \textbf{88} (1984), 279--291;
Journal of Algebra \textbf{93} (1985), 310--323.

\bibitem{Ku2}
 Burkhard K\"ulshammer,
 {\em Group-theoretical descriptions of ring theoretical
    invariants of group algebras},
 Progress in Mathematics, \textbf{95} (1991), 425--441.

\bibitem{Liu2008}
Yuming Liu, {\em Summands of stable equivalences of Morita type}.
Communications in Algebra \textbf{36}(10) (2008), 3778-3782.

\bibitem{LZZ}
Yuming Liu, Guodong Zhou and Alexander Zimmermann, {\em Higman
ideal, stable Hochschild homology and Auslander-Reiten
conjecture}, preprint (2008).

\bibitem{Loday}
Jean-Louis Loday, {\sc Cyclic Homology}, Springer Verlag, Berlin 1992.

\bibitem{Murray1}
John Murray, {\em On a certain ideal of K\"ulshammer in the centre of a group algebra},
Archiv der Mathematik  {\bf 77} (2001) 373-377.

\bibitem{Murray2}
John Murray, {\em Squares in the centre of the group algebra of a symmetric group},
Bulletin of the London Mathematical Society {\bf 34} (2002) 155-164.

\bibitem{Nakayama}
Tadasi Nakayama, {\em On Frobeniusean algebras I,\,II},
Annals of Mathematics {\bf 40} (1939) 611-633;
Annals of Mathematics {\bf 42} (1941) 1-21.

\bibitem{Pogorzaly}
Zygmunt Pogorza\l y, {\em Algebras stably equivalent to
self-injective special biserial algebras},  Communications in Algebra
\textbf{22}  (1994),  no. 4, 1127-1160

\bibitem{Quillen}
Daniel Quillen, {\em Cyclic cohomology and algebra extensions}, $K$-Theory
{\bf 3} (1989) 205-246.

\bibitem{MO}
Irving Reiner, {\sc Maximal Orders}, Academic Press London, New York, San Francisco 1975

\bibitem{Reynolds}
William F. Reynolds, {\em Sections and ideals of centres of group algebras},
Journal of Algebra {\bf 20} (1972) 176-181.

\bibitem{Rickard1989} Jeremy Rickard,
\emph{ Derived categories and stable equivalence}.
Journal of Pure and Applied Algebra  \textbf{61}  (1989),  no. 3, 303--317.


\bibitem{Ri1}
Jeremy Rickard, {\em Derived equivalences as derived functors},
Journal of the London Mathematical Society {\bf 43} (1991) 37-48.

\bibitem{Stasheff}
Jim Stasheff, {\em The intrinsic bracket on the deformation complex
of an associative algebra}, Journal of Pure and Applied Algebra {\bf 89}
(1993) 231-235.

\bibitem{Verdier}
Jean-Louis Verdier, {\sc Des cat\'egories d\'eriv\'ees des cat\'egories ab\'eliennes}, Ast\'erisque {\bf 239} Soci\'et\'e math\'ematique de France 1996.

\bibitem{Xi2008} Chang-Chang Xi,
{\em Stable equivalences of adjoint type}. Forum Mathematicum \textbf{20}
(2008), 81-97.

\bibitem{Yamagata}
Kunio Yamagata, {\em Frobenius algebras}, in Handbook of algebra Vol.1, 1996, pages 841-887, North-Holland Amsterdam.

\bibitem{Modulestructure}
Guodong Zhou and Alexander Zimmermann, {\em Classifying tame blocks and
related algebras up to stable equivalences of Morita type}, preprint 2010

\bibitem{polygrowth}
Guodong Zhou and Alexander Zimmermann, {\em Auslander-Reiten conjecture
for symmetric algebras of polynomial growth}, preprint 2010

\bibitem{Mexico} Alexander Zimmermann,
{\em Derived Equivalences of Orders}, in Proceedings of the ICRA VII,
 Mexico, eds: Bautista, Martinez, de la Pe\~{n}a.
 Canadian Mathematical Society Conference Proceedings \textbf{18}, (1996)
 721-749.

\bibitem{rogquest}
Alexander Zimmermann,
{\em Tilted symmetric orders are symmetric orders},
Archiv der Mathematik {\bf 73 } (1999) 15-17.

\bibitem{Kuelsquest}
Alexander Zimmermann, {\em Invariance of generalized Reynolds ideals
under derived equivalences.}
Mathematical Proceedings of the Royal Irish Academy {\bf 107}A (1)
(2007) 1-9.

\bibitem{Gerstenhaber}
Alexander Zimmermann,
{\em Fine Hochschild invariants of derived categories for symmetric
algebras},
Journal of Algebra {\bf 308} (2007) 350-367

\bibitem{TAHochschild}
Alexander Zimmermann,
{\em Hochschild
homology invariants of K\"ulshammer type of derived categories},
to appear in Communications in Algebra.

\end{thebibliography}
\end{document}